\input amstex 
\font\we=cmb10 at 14.4truept
\font\li=cmb10 at 12truept
\documentstyle{amsppt} 
\catcode`\@=11
\redefine\logo@{}
\catcode`\@=13
\hsize=6.5truein
\vsize=9.5truein
\topmatter 
\title\nofrills {\we Non-Abelian  Zeta Functions For Function Fields} 
\endtitle 
\author {\li Lin WENG} 
\endauthor 
\affil {\bf Graduate School of Mathematics, Kyushu University, Japan}\endaffil
\NoRunningHeads 
\endtopmatter 
\TagsOnRight
\vskip 0.45cm
\noindent
{\bf Abstract.} In this paper we initiate a geometrically oriented 
construction of
non-abelian zeta functions for curves defined over finite fields.
More precisely, we first introduce new yet genuine non-abelian
zeta functions for curves  defined over finite fields, by  a
\lq weighted count\rq$\ $on rational points over the corresponding
moduli spaces of semi-stable vector bundles using moduli 
interpretation of these points. Then we define non-abelian $L$-functions 
for  curves over finite fields
using integrations of Eisenstein series associated to 
$L^2$-automorphic forms over certain generalized  
moduli spaces. 
\vskip 1.0cm
In this paper we initiate a geometrically oriented construction of
non-abelian zeta functions for curves defined over finite fields. 
It consists of two chapters.

More precisely, in Chapter I, we first introduce new yet genuine non-abelian
zeta functions for curves  defined over finite fields. This is achieved by  a
\lq weighted count\rq$\ $on rational points over the corresponding
moduli spaces of semi-stable vector bundles using moduli 
interpretation of these points. We justify our construction by 
establishing basic properties for these new zetas such as  functional 
equation and rationality, and show that if only line bundles 
are involved,  our newly defined zetas coincide with  Artin's Zeta.
All this, in particular, the rationality, then leads naturally to  
our definition of (global) non-abelian zeta functions 
(for curves defined over number fields), which themselves are
justified by a convergence result. 
We end this chapter with a detailed study on rank two 
non-abelian zeta functions for genus two curves, based on
what we call infinitesimal structures of Brill-Noether loci (and 
Weierstrass points).

In Chapter II, we begin with a similar construction for the field of rationals
to motivate what follows. In particular, we show that  there is an 
intrinsic relation 
between our non-abelian zeta functions and Eisenstein series. 
Due to this, instead of introducing general non-abelian $L$-functions
for curves defined over finite fields with more general test functions 
(as what Tate did in his Thesis for abelian $L$-functions), 
we then define non-abelian $L$-functions for  curves over finite fields
as integrations of Eisenstein series associated to 
$L^2$-automorphic forms over certain generalized  
moduli spaces. Here geometric truncations play a key role. 
Basic properties for these non-abelian $L$-functions, such as 
meromorphic continuation, functional equations and singularities,
are established as well, based on the theory of Eisenstain series of 
Langlands and Morris.  We end this chapter by establishing 
a closed formula for what we call the abelian parts of
non-abelian $L$-functions associated with Eisenstein series for cusp 
forms, via the
Rankin-Selberg method, motivated by a formula of Arthur and Langlands.

This work is an integrated part of our vast yet still under developing 
Program for Geometric Arithmetic [We1], and 
is motivated by our new non-abelian $L$-functions 
for number fields [We2] in connection with non-abelian arithmetic 
aspects of global fields.

\vskip 0.50cm
\centerline {\li Chapter I. Non-Abelian Zeta Functions}
\vskip 0.30cm
This consists of two aspects: construction and justification.
For the construction,  we first introduce a new type of zeta functions 
for curves defined over finite fields using the corresponding moduli 
spaces of semi-stable vector bundles. We show that these new zeta
functions are indeed rational and satisfy certain functional equation, based
on vanishing theorem, (duality, Riemann-Roch theorem) for cohomologies 
of semi-stable vector bundles. Based on this, in particular, the rationality,  
we then introduce global non-abelian zeta functions for curves
defined over number fields, via the  Euler product formalism. 
Moreover, we establish a convergence result
for our Euler products using the Clifford Lemma, 
an ugly yet quite explicit formula for local non-abelian
zeta functions,  a  result of (Harder-Narasimhan) Siegel about quadratic
forms, and Weil's theorem on Riemann Hypothesis for Artin zeta functions.

As for the justification, surely, we check that when only line bundles are
involved, (so  moduli spaces of semi-stable bundles are nothing but the 
standard  Picard groups), our (new) zeta functions, global and local, 
coincide with the classical Artin zeta functions for curves defined over finite
fields  and Hasse-Weil zeta functions for curves defined over number 
fields respectively. 
Moreover, as concrete examples, 
we compute  rank two zeta functions for genus two curves 
by studying Weierstrass points and non-abelian Brill-Noether loci 
in terms of what we call their infinitesimal structures. 
\vskip 0.45cm
\centerline {\li I.1. Local Non-Abelian Zeta Functions for Curves}
\vskip 0.30cm 
In this section, we introduce our non-abelian
zeta functions for curves defined over finite fields. Basic properties 
for these non-abelian zeta functions, such as meromorphic
extensions, rationality and functional equations, are established.
\vskip 0.30cm
\centerline {\bf 1.1. Moduli Spaces of Semi-Stable Bundles}
\vskip 0.30cm
\noindent
{\bf 1.1.1. Semi-Stable Bundles.} Let $C$ be a regular, reduced and 
irreducible projective curve defined over an algebraically closed field
$\bar k$. Then according to Mumford [Mu], a vector bundle $V$ on 
$C$ is called semi-stable (resp. stable) if for any proper subbundle 
$V'$ of $V$,  $$\mu(V'):={{d(V')}\over {r(V')}}\leq \text{
(resp.}< \text{)}{{d(V)}\over {r(V)}}=:\mu(V).$$ Here $d$ denotes 
the degree and $r$ denotes the
rank.
\vskip .30cm
\noindent
{\bf Proposition.} {\it Let $V$ be a vector bundle over $C$. Then} 

\noindent
(a) ([HN]) {\it there exists a unique
filtration of subbundles of $V$, the  Harder-Narasimhan filtration of $V$, 
$$\{0\}=V_0\subset V_1\subset V_2\subset\dots\subset V_{s-1}\subset V_s=V$$
such that all $V_i/V_{i-1}$ are semi-stable and for $1\leq i\leq s-1$, 
$\mu(V_i/V_{i-1})>\mu(V_{i+1}/V_i);$}

\noindent
(b) (see e.g. [Se]) {\it if moreover $V$ is semi-stable,
there exists a
filtration of subbundles of $V$, a Jordan-H\"older filtration of $V$,
$$\{0\}=V^{t+1}\subset V^t\subset\dots\subset V^{1}\subset V^0=V$$
such that for all $0\leq i\leq t$, $V^i/V^{i+1}$ is stable and 
$\mu(V^i/V^{i+1})
=\mu (V)$. Moreover,
the associated graded bundle $\text{Gr}(V):=\oplus_{i=0}^tV^i/V^{i+1}$, the 
(Jordan-H\"older) graded bundle of
$V$, is determined uniquely by $V$.}
\vskip 0.30cm
\noindent
{\bf 1.1.2. Moduli Space of Stable Bundles.} Following Seshadri, 
two semi-stable 
vector bundles $V$
and $W$ are called $S$-equivalent, if their associated Jordan-H\"older graded 
bundles are isomorphic, i.e.,
$\text{Gr}(V)\simeq \text{Gr}(W)$. Applying Mumford's general result on 
geometric invariant theory, Narasimhan and Seshadri proved the following 
\vskip 0.30cm
\noindent
{\bf Theorem.} (See e.g. [NS] and [Se]) {\it Let $C$ be a regular, 
reduced, irreducible projective 
curve of genus $g\geq 2$ defined over an algebraically closed field. 
Then over the set ${\Cal M}_{C,r}(d)$ (resp. ${\Cal M}_{C,r}(L)$) of 
$S$-equivalence classes of rank $r$ and degree $d$ (resp. rank $r$ and
determinant $L$) semi-stable vector bundles over $C$, there is a 
natural normal, 
projective $(r^2(g-1)+1)$-dimensional (resp. $(r^2-1)(g-1)$-dimensional) 
algebraic variety structure.}
\vskip 0.30cm
\noindent
{\it Remark.} In this paper, we always assume that the genus of $g$ is 
at least 2. 
For  elliptic curves, whose associated moduli spaces are very special, 
please see [We3].
\vskip 0.30cm
\noindent
{\bf 1.1.3. Rational Points.} Now  assume that $C$ is defined over a finite 
field $k$. It makes sense to talk about
$k$-rational bundles over $C$, i.e., bundles which are defined over $k$. 
Moreover, from geometric invariant theory, projective varieties
${\Cal M}_{C,r}(d)$ are defined over a certain finite extension of $k$; 
and if $L$ itself is defined over $k$, the same holds for ${\Cal M}_{C,r}(L)$. 
Thus it makes sense to talk about $k$-rational points of these moduli 
spaces too. 
The relation between these two types of rationality
is given by  Harder-Narasimhan based on a discussion about Brauer groups: 
\vskip 0.30cm
\noindent
{\bf Proposition.} ([HN]) {\it Let $C$ be a regular, reduced, irreducible 
projective curve of genus $g\geq 2$ defined over a finite field $k$. 
Then there exists a finite field ${\Bbb F}_q$ such that for all
$d$ (resp. all $k$-rational line bundles $L$), the subset of 
${\Bbb F}_q$-rational points of ${\Cal M}_{C,r}(d)$ 
(resp. ${\Cal M}_{C,r}(L)$) consists exactly of all $S$-equivalence 
classes of
${\Bbb F}_q$-rational bundles in ${\Cal M}_{C,r}(d)$ (resp. ${\Cal
M}_{C,r}(L)$).}
\vskip 0.30cm
From now on, without loss of generality, we always assume that 
the finite fields ${\Bbb F}_q$ (with $q$ elements) satisfy the 
property stated in the Proposition. Also for simplicity, we 
write ${\Cal M}_{C,r}(d)$ (resp. ${\Cal M}_{C,r}(L)$) for 
${\Cal M}_{C,r}(d)({\Bbb F}_q)$ (resp. ${\Cal M}_{C,r}(L)({\Bbb F}_q)$), 
the subset of ${\Bbb F}_q$-rational points, and call them moduli spaces 
by an abuse of notations. Clearly these sets are all finite.
\vskip 0.45cm
\centerline {\bf 1.2. Local Non-Abelian Zeta Functions}
\vskip 0.30cm
\noindent
{\bf 1.2.1. Definition.} Let $C$ be a  regular, reduced, irreducible 
projective 
curve of genus $g\geq
2$ defined over the finite field ${\Bbb F}_q$ with $q$ elements. Define 
the {\it rank $r$
non-abelian zeta function  $\zeta_{C,r,{\Bbb F}_q}(s)$ of $C$}
by setting
$$\zeta_{C,r,{\Bbb F}_q}(s):=\sum_{V\in [V]\in {\Cal M}_{C,r}(d), 
d\geq 0}{{q^{h^0(C,V)}-1}\over
{\#\text{Aut}(V)}}\cdot (q^{-s})^{d(V)},\qquad \text{Re}(s)>1.$$
\vskip 0.30cm
\noindent
{\bf Proposition.} {\it With the same notation as above, 
$\zeta_{C,1,{\Bbb F}_q}(s)$ is nothing but
the classical Artin zeta function $\zeta_C(s)$ for curve $C$. 
That is to say,
$$\zeta_{C,1,{\Bbb F}_q}(s)=\sum_{D\geq 0}{1\over {N(D)^{s}}}=:\zeta_C(s)
\qquad \text{Re}(s)>1.$$ Here $D$ runs
over all effective divisors of $C$, and $N(D):=q^{d(D)}$ with
$d(\Sigma_Pn_PP):=\Sigma_Pn_Pd(P)$.}
\vskip 0.30cm
\noindent
{\it Proof.} By definition, the classical Artin zeta function ([A], [Mo]) 
for $C$ is given by
$$\zeta_C(s):=\sum_{D\geq 0}{1\over {N(D)^{s}}}.$$ 
 Thus by first grouping effective divisors 
according to their rational equivalence classes ${\Cal D}$, then taking the 
sum on effective divisors in the same class,
we obtain $$\zeta_C(s)=\sum_{\Cal D}\sum_{D\in {\Cal D},D \geq 0}
{1\over {N(D)^{s}}}.$$ Clearly, $$\sum_{D\in {\Cal D}, D \geq 0}
{1\over {N(D)^{s}}}={{q^{h^0(C,{\Cal D})}-1}\over
{q-1}}\cdot (q^{-s})^{d({\Cal D})}.$$ Therefore, 
$$\zeta_C(s)=\sum_{L\in \text{
Pic}^d(C), d\geq 0}{{q^{h^0(C,L)}-1}\over {\#\text{Aut}(L)}}
\cdot (q^{-s})^{d(L)}$$ due to the fact that
$\text{Aut}(L)\simeq {\Bbb F}_q^*$.
\vskip 0.30cm
\noindent
{\it Remark.} Before going further, let us explain the notation $V\in [V]$
 appeared in the summation in detail. By $\sum_{V\in [V]}$, we mean that 
the sum is taking over  
all (isomorphism classes of) rational vector bundles $V$ in $[V]$. 
From  Prop. (b) in 1.1.1,  for each fixed $[V]$, there are only finitely 
many
 terms involved. On the other hand, we may instead use only a single element 
$V$ for each class $[V]$, say, one with maximal automorphism group 
(as used in the proof of the projectivity of  moduli spaces). However,
while interesting, such a change yields  quite different functions. 
(See e.g. [We1].)
Our decision to use  all rational elements in $[V]$
is motivated by an adelic consideration, in particular, by
Harder-Narasimhan's understanding of Siegel's formula. 
\vskip 0.30cm
\noindent
{\bf 1.2.2. Convergence and Rationality.} At this point, we must show that for 
general $r$, the infinite
summation in the definition of our non-abelian zeta function 
$\zeta_{C,r,{\Bbb F}_q}(s)$ converges
when $\text{Re}(s)>1$. For this, let us start with  the following simple
vanishing result for semi-stable vector bundles.

\noindent
{\bf Lemma 1.} {\it Let $V$ be a  rank $r$ semi-stable vector bundle of 
degree $d$
on $C$. Then 

\noindent
(a) if $d\geq r(2g-2)+1$, $h^1(C,V)=0$;

\noindent
(b) if $d<0$, $h^0(C,V)=0$.}
\vskip 0.30cm
\noindent
{\it Proof.} This is a direct consequence of the fact that if $V$ and $W$ 
are semi-stable vector bundles with $\mu(V)>\mu(W)$, 
then $H^0(C,\text{Hom}(V,W))=\{0\}$. 
\vskip 0.30cm
Thus, from definition, 
$$\eqalign{\zeta_{C,r,{\Bbb F}_q}(s)=&
\sum_{V\in [V]\in {\Cal M}_{C,r}(d), 0\leq d\leq r(2g-2)
}{{q^{h^0(C,V)}-1}\over {\#\text{Aut}(V)}}\cdot (q^{-s})^{d(V)}\cr
&\qquad+\sum_{V\in [V]\in {\Cal
M}_{C,r}(d), d\geq r(2g-2)+1}{{q^{d(V)-r(g-1)}-1}\over 
{\#\text{Aut}(V)}}\cdot (q^{-s})^{d(V)}.\cr}$$
Clearly only finitely many terms appear in the first summation, 
so it suffices to show that when
$\text{Re}(s)>1$, the second term converges. For this purpose, 
we introduce what we call the Harder-Narasimhan numbers 
$$\beta_{C,r,{\Bbb F}_q}(d):=\sum_{V\in [V]\in {\Cal M}_{C,r}(d)}
{1\over {\#\text{Aut}(V)}}.$$ 

\noindent
{\bf Lemma 2.} {\it With the same notation as above, for all $n\in {\Bbb Z}$,
$$\beta_{C,r,{\Bbb F}_q}(d+rn)=\beta_{C,r,{\Bbb F}_q}(d).$$}
{\it Proof.} This comes from the following two facts: 

\noindent
(1) there is a degree one ${\Bbb F}_q$-rational
line bundle $A$ on $C$; and 

\noindent
(2) $\text{Aut}(V)\simeq
\text{Aut}(V\otimes A^{\otimes n})$ and $d(V\otimes A^{\otimes n})=d(V)+rn$.
\vskip 0.30cm
Therefore,  the second summation becomes
$$\eqalign{~&\sum_{i=1}^{r}\beta_{C,r,{\Bbb F}_q}(i)
\sum_{n=2g-2}^\infty\Big(q^{nr+i-r(g-1)}-1\Big)\cdot
(q^{-s})^{nr+i}\cr
=&\sum_{i=1}^r\beta_{C,r,{\Bbb F}_q}(i) \cdot (q^{-s})^i\cdot
\Big(q^{i-r(g-1)}\cdot {{q^{(1-s)\cdot r(2g-2)}}\over {1-q^{(1-s)\cdot r}}}
-{{q^{(-s)\cdot r(2g-2)}}\over {1-q^{(-s)\cdot r}}}\Big),\cr}$$ 
provided that $|q^{-s}|<1$. Thus we have proved the following
\vskip 0.30cm
\noindent
{\bf Proposition.} {\it The non-abelian zeta function 
$\zeta_{C,r,{\Bbb F}_q}(s)$ is well-defined
for $\text{Re}(s)>1$, and admits a meromorphic extension 
to the whole complex $s$-plane.}
\vskip 0.30cm
Moreover, if we  set
$t:=q^{-s}$ and introduce the non-abelian $Z$-function of $C$ by 
$$Z_{C,r,{\Bbb F}_q}(t):=
\sum_{V\in [V]\in {\Cal M}_{C,r}(d),d\geq
0}{{q^{h^0(C,V)}-1}\over {\#\text{Aut}(V)}}\cdot t^{d(V)}, \qquad |t|<1.$$ 
Then the above calculation implies that
$$Z_{C,r,{\Bbb F}_q}(t)
=\sum_{d=0}^{r(2g-2)}\Big(\sum_{V\in [V]\in {\Cal
M}_{C,r}(d)}{{q^{h^0(C,V)}-1}\over {\#\text{Aut}(V)}}\Big)\cdot t^{d}
+\sum_{i=1}^r\beta_{C,r,{\Bbb
F}_q}(i)\cdot \Big({{q^{r(g-1)+i}}\over {1-q^rt^r}}-{{1}
\over {1-t^r}}\Big)\cdot t^{r(2g-2)+i}.$$
Therefore, there exists a polynomial 
$P_{C,r,{\Bbb F}_q}(s)\in {\Bbb Q}[t]$ such that
$$Z_{C,r,{\Bbb F}_q}(t)={{P_{C,r,{\Bbb F}_q}(t)}\over {(1-t^r)(1-q^rt^r)}}.$$ 
In this way, we have established the following
\vskip 0.30cm
\noindent 
{\bf Rationality.} {\it Let $C$ be a regular, reduced irreducible 
projective curve defined over 
${\Bbb F}_q$ with $Z_{C,r,{\Bbb F}_q}(t)$ the rank $r$ non-abelian 
$Z$-function. Then, there exists a polynomial
$P_{C,r,{\Bbb F}_q}(s)\in {\Bbb Q}[t]$ such that 
$$Z_{C,r,{\Bbb F}_q}(t)={{P_{C,r,{\Bbb F}_q}(t)}\over {(1-t^r)(1-q^rt^r)}}.$$}

\noindent
{\bf 1.2.3. Functional Equation.} To understand
$P_{C,r,{\Bbb F}_q}(s)$ better, as well as for  theoretical purpose, 
we next study  functional
equation for  rank $r$ zeta functions. Let us 
introduce the rank $r$ non-abelian $\xi$-function
$\xi_{C,r,{\Bbb F}_q}(s)$ by setting
$$\xi_{C,r,{\Bbb F}_q}(s):=\zeta_{C,r,{\Bbb F}_q}(s)\cdot (q^{s})^{r(g-1)}.$$ 
That is to say,
$$\xi_{C,r,{\Bbb F}_q}(s)=\sum_{V\in [V]\in {\Cal M}_{C,r}(d),d\geq 0}
{{q^{h^0(C,V)}-1}\over {\#\text{Aut}(V)}}\cdot (q^{-s})^{\chi(C,V)}, 
\qquad \text{Re}(s)>1,$$
where $\chi(C,V)$ denotes the Euler-Poincar\'e characteristic of $V$.
\vskip 0.30cm
\noindent
{\bf Functional Equation.} {\it Let $C$ be a regular, reduced irreducible 
projective curve defined over
${\Bbb F}_q$ with $\xi_{C,r,{\Bbb F}_q}(s)$ its associated rank $r$ 
non-abelian  $\xi$-function. Then,
$$\xi_{C,r,{\Bbb F}_q}(s)=\xi_{C,r,{\Bbb F}_q}(1-s).$$}

Before proving  the functional equation, we give the following
\vskip 0.30cm
\noindent
{\bf Corollary.} {\it With the same notation as above,

\noindent
(a) $P_{C,r,{\Bbb F}_q}(t)\in {\Bbb Q}[t]$ is a degree $2rg$ polynomial;

\noindent
(b) Denote all reciprocal roots of $P_{C,r,{\Bbb F}_q}(t)$ by 
$\omega_{C,r,{\Bbb F}_q}(i),
i=1,\dots, 2rg$. Then after a suitable rearrangement,
$$\omega_{C,r,{\Bbb F}_q}(i)\cdot \omega_{C,r,{\Bbb F}_q}(2rg-i)=q,
\qquad i=1,\dots,rg;$$

\noindent
(c) For each $m\in {\Bbb Z}_{\geq 1}$, there exists a rational 
number $N_{C,r,{\Bbb F}_q}(m)$ such that
$$Z_{r,C,{\Bbb F}_q}(t)=P_{C,r,{\Bbb F}_q}(0)\cdot
\exp\Big(\sum_{m=1}^\infty N_{C,r,{\Bbb
F}_q}(m){{t^m}\over m}\Big).$$ Moreover, $$N_{C,r,{\Bbb
F}_q}(m)=\cases r(1+q^m)-\sum_{i=1}^{2rg}
\omega_{C,r,{\Bbb F}_q}(i)^m,& r\ |m;\\
-\sum_{i=1}^{2rg}\omega_{C,r,{\Bbb F}_q}(i)^m,& 
r\not| m;\endcases$$

\noindent
(d) For any  $a\in {\Bbb Z}_{>0}$, denote by $\zeta_{a}$ a primitive 
$a$-th root of unity and set $T=t^a$.  Then
$$\prod_{i=1}^aZ_{C,r}(\zeta_{a}^it)=(P_{C,r,{\Bbb F}_q}(0))^a\cdot
\exp\Big(\sum_{m=1}^\infty N_{r,C,{\Bbb F}_q}(ma){{T^m}\over m}\Big).$$}

\noindent
{\it Proof.} (a) and (b) are direct consequences of the functional equation, 
while (c) and (d) are
direct consequences of (a), (b) and the following well-known relations
$$\sum_{i=1}^a(\zeta_{a}^i)^m=\cases a,& a\ |m,\\
0, & a\not|m.\endcases$$

\noindent
{\bf 1.2.4. Proof of Functional Equation.} To understand the structure of 
the functional equation
explicitly, we decompose the non-abelian $\xi$-function for curves.
For this purpose, first recall that the canonical line bundle $K_C$ of 
$C$ is defined over ${\Bbb F}_q$. Thus, for all
$n\in {\Bbb Z}$, we obtain the following natural ${\Bbb F}_q$-rational 
isomorphisms:
$$\matrix {\Cal M}_r(L)&\to& {\Cal M}_r(L\otimes K_C^{\otimes nr});&\qquad&
{\Cal M}_r(L)&\to& {\Cal M}_r(L^{\otimes -1}\otimes K_C^{\otimes nr})\\
[V]&\mapsto&[V\otimes K_C^{\otimes n}];&\qquad&[V]&\mapsto& 
[V^\vee\otimes K_C^{\otimes n}],\endmatrix$$  where  $V^\vee$ 
denotes the dual of $V$. Next, introduce the union  
$${\Cal M}_{C,r}^L:=\cup_{n\in {\Bbb Z}}\Big({\Cal M}_r(L\otimes
K_C^{\otimes nr})\cup {\Cal M}_r(L^{\otimes -1}\otimes K_C^{\otimes
nr})\Big).$$ With this, clearly, we may and indeed always assume that
$$0\leq d(L)\leq r(g-1).$$

Furthermore, introduce the partial non-abelian zeta function 
$\xi_{C,r,{\Bbb F}_q}^L(s)$ by setting
$$\xi_{C,r,{\Bbb F}_q}^L(s):=\sum_{V\in [V]\in 
{\Cal M}_{C,r}^L}{{q^{h^0(C,V)}-1}\over {\#\text{
Aut}(V)}} \cdot\big(q^{-s}\big)^{\chi(C,V)},\qquad \text{Re}(s)>1.$$ 
Clearly, then
$$\xi_{C,r,{\Bbb F}_q}(s)=\sum_L\xi_{C,r,{\Bbb F}_q}^L(s)$$ where $L$ 
runs over all line bundles
appeared in the following (disjoint) union
$$\cup_{d\in {\Bbb Z}}{\Cal M}_{C,r}(d)=\cup_{L,0\leq d(L)
\leq r(g-1)}{\Cal M}_{C,r}^L.$$
 Here we reminder the reader that the vanishing result 
of Lemma 1.2.2.1 has been used.
 
Therefore, to establish the functional equation for $\xi_{C,r,{\Bbb F}_q}(s)$, 
it suffices to show that
$$\xi_{C,r,{\Bbb F}_q}^L(s)=\xi_{C,r,{\Bbb F}_q}^L(1-s).$$
For this, we have the following

\noindent
{\bf Theorem.} {\it  For $\text{Re}(s)>1$,  
$$\eqalign{~&\xi_{C,r,{\Bbb F}_q}^L(s)\cr
=&{1\over
2}\sum_{V\in [V]\in {\Cal M}_{C,r}^L;0\leq d(V)\leq
r(2g-2)}{{q^{h^0(C,V)}}\over {\#\text{
Aut}(V)}}\cdot\Big[(q^{-s})^{\chi(C,V)}+
(q^{s-1})^{\chi(C,V)}\Big]\cr
&+\Big[{{q^{(1-s)\cdot(d(L)-r(g-1))}}\over
{q^{(s-1)\cdot r(2g-2)}-1}}
+{{q^{s\cdot(d(L)-r(g-1))}}\over
{q^{(-s)\cdot r(2g-2)}-1}}+{{q^{(s-1)\cdot(d(L)-r(g-1))}}\over
{q^{(s-1)\cdot r(2g-2)}-1}}
+{{q^{(-s)\cdot(d(L)-r(g-1))}}\over
{q^{(-s)\cdot r(2g-2)}-1}}\Big]\cdot \beta_{C,r,{\Bbb F}_q}(L).\cr}\eqno(*)$$ 
Here
$\beta_{C,r,{\Bbb F}_q}(L):=\sum_{V\in [V]\in {\Cal M}_{C,r}(L)}
{1\over{\#\text{
Aut}(E)}}$ denotes the Harder-Narasimhan number. In particular, 

\noindent
(a) $\xi_{C,r,{\Bbb F}_q}^L(s)$ satisfies the functional equation
$\qquad\xi_{C,r,{\Bbb F}_q}^L(s)=\xi_{C,r,{\Bbb F}_q}^L(1-s);$

\noindent
(b) the Harder-Narasimhan number $\beta_{C,r,{\Bbb F}_q}(L)$ is given 
by the leading
term of the  singularities of $\xi_{C,r,{\Bbb F}_q}^L(s)$  at $s=0$ and $s=1$.}
\vskip 0.30cm
\noindent
{\it Proof.} It suffices to prove (*). For this, set
$$I(s)=\sum_{V\in [V]\in {\Cal M}_{C,r}^L;0\leq d(V)\leq
r(2g-2)}{{q^{h^0(C,V)}}\over {\#\text{
Aut}(V)}}\cdot (q^{-s})^{\chi(C,V)}$$ and
$$II(s)=\sum_{V\in [V]\in {\Cal M}_{C,r}^L;=d(V)>
r(2g-2)}{{q^{h^0(C,V)}}\over {\#\text{
Aut}(V)}}\cdot (q^{-s})^{\chi(C,V)}-\sum_{V\in [V]\in 
{\Cal M}_{C,r}^L;d(V)\geq 0}{{1}\over {\#\text{
Aut}(V)}}\cdot (q^{-s})^{\chi(C,V)}.$$ Thus, 
$$\xi_{C,r,{\Bbb F}_q}^L(s)=I(s)+II(s).$$  So it suffices to show the following

\noindent
{\bf Lemma.} {\it With the same notation as above, 

\noindent
(a) $I(s)={1\over
2}\sum_{V\in [V]\in {\Cal M}_{C,r}^L;0\leq d(V)\leq
r(2g-2)}{{q^{h^0(C,V)}}\over {\#\text{
Aut}(V)}}\cdot\Big[(q^{-s})^{\chi(C,V)}+
(q^{s-1})^{\chi(C,V)}\Big];$ and

\noindent
(b) $$II(s)=\Big[{{q^{(1-s)\cdot(d(L)-r(g-1))}}\over
{q^{(s-1)\cdot r(2g-2)}-1}}
+{{q^{s\cdot(d(L)-r(g-1))}}\over
{q^{(-s)\cdot r(2g-2)}-1}}+{{q^{(s-1)\cdot(d(L)-r(g-1))}}\over
{q^{(s-1)\cdot r(2g-2)}-1}}
+{{q^{(-s)\cdot(d(L)-r(g-1))}}\over
{q^{(-s)\cdot r(2g-2)}-1}}\Big]\cdot \beta_{C,r,{\Bbb F}_q}(L).$$}

\noindent
{\it Proof.} (a) comes from  Riemann-Roch theorem and  Serre duality. Indeed,
$$\eqalign{I(s)
=&{1\over
2}\Big(\sum_{V\in [V]\in {\Cal M}_{C,r}^L;0\leq d(E)\leq
r(2g-2)}{{q^{h^0(C,V)}}\over {\#\text{
Aut}(V)}}\cdot(q^{-s})^{\chi(C,V)}\cr
&\qquad+\sum_{V^\vee\otimes
K_C\in {\Cal M}_{C,r}^L;0\leq d(V^\vee\otimes K_C)\leq
r(2g-2)}{{q^{h^0(C,V^\vee\otimes K_C)}}\over {\#\text{
Aut}(V^\vee\otimes K_C)}}\cdot(q^{-s})^{\chi(C,V^\vee\otimes
K_C)}\Big)\cr
=&{1\over
2}\sum_{V\in [V]\in {\Cal M}_{C,r}^L;0\leq d(E)\leq
r(2g-2)}\Big[{{q^{h^0(C,V)}}\over {\#\text{
Aut}(V)}}\cdot(q^{-s})^{\chi(C,V)}+{{q^{h^1(C,V^\vee\otimes K_C)}}
\over {\#\text{
Aut}(V^\vee\otimes K_C)}}\cdot(q^{1-s})^{\chi(C,V^\vee\otimes
K_C)}\Big]\cr
=&{1\over
2}\sum_{V\in [V]\in {\Cal M}_{C,r}^L;0\leq d(V)\leq
r(2g-2)}{{q^{h^0(C,V)}}\over {\#\text{
Aut}(V)}}\cdot\Big[(q^{-s})^{\chi(C,V)}+
(q^{s-1})^{\chi(C,V)}\Big].\cr}$$ 

As for (b), clearly by the vanishing result,
$$\eqalign{T_{r,L}(s)=&\sum_{V\in [V]\in {\Cal M}_{C,r}^L;d(E)>
r(2g-2)}{1\over {\#\text{
Aut}(E)}}\cdot(q^{1-s})^{\chi(C,V)}-\sum_{V\in [V]\in 
{\Cal M}_{C,r}^L;d(E)\geq 0}
{1\over {\#\text{Aut}(E)}}\cdot(q^{-s})^{\chi(C,V)}\cr
=&\Big(\sum_{V\in [V]\in {\Cal M}_{C,r}(L\otimes K_C^{\otimes rn});
d(L)+rn(2g-2)> r(2g-2)}
{1\over {\#\text{Aut}(V)}}\cdot(q^{1-s})^{\chi(C,E)}\cr
&\qquad-\sum_{V\in [V]\in {\Cal M}_{C,r}(L^{-1}\otimes K_C^{\otimes rn});
-d(L)+rn(2g-2)\geq 0}
{1\over {\#\text{Aut}(V)}}\cdot(q^{-s})^{\chi(C,E)}\Big)\cr
&+\Big(\sum_{V\in [V]\in {\Cal M}_{C,r}(L^{-1}\otimes K_C^{\otimes rn});
-d(L)+rn(2g-2)> r(2g-2)}
{1\over {\#\text{Aut}(V)}}\cdot(q^{1-s})^{\chi(C,E)}\cr
&\qquad-\sum_{V\in [V]\in {\Cal M}_{C,r}(L\otimes K_C^{\otimes rn});
d(L)+rn(2g-2)>0}
{1\over {\#\text{Aut}(V)}}\cdot(q^{-s})^{\chi(C,E)}\Big).\cr}$$
But $\chi(C,V)$ depends only on $d(V)$. Thus, accordingly,
$$\eqalign{II(s)
=&\Big[\Big(\sum_{n=1}^\infty (q^{1-s})^{d(L)+nr(2g-2)-r(g-1)}
-\sum_{n=1}^\infty(q^{-s})^{-d(L)+nr(2g-2)-r(g-1)}\Big)\cr 
&+\Big(\sum_{n=2}^\infty(q^{1-s})^{-d(L)+nr(2g-2)-r(g-1)}
-\sum_{n=0}^\infty(q^{-s})^{d(L)+nr(2g-2)-r(g-1)}\Big)\Big]\cdot
\beta_{C,r}(L)\cr 
=&\Big[{{q^{(1-s)\cdot(d(L)-r(g-1))}}\over
{q^{(s-1)\cdot r(2g-2)}-1}}
+{{q^{s\cdot(d(L)-r(g-1))}}\over
{q^{(-s)\cdot r(2g-2)}-1}}+{{q^{(s-1)\cdot(d(L)-r(g-1))}}\over
{q^{(s-1)\cdot r(2g-2)}-1}}
+{{q^{(-s)\cdot(d(L)-r(g-1))}}\over
{q^{(-s)\cdot r(2g-2)}-1}}\Big]\cdot \beta_{C,r,{\Bbb F}_q}(L).\cr}$$ 
This completes the proof of the lemma, and hence the Theorem and the
Functional Equation for rank $r$ zeta functions.
\eject
\vskip 0.45cm
\centerline {\li I.2. Global Non-Abelian Zeta Functions for Curves}
\vskip 0.30cm
In this section, we introduce new non-abelian zeta functions for 
curves defined over number fields 
via the Euler product formalism, based on our study of
non-abelian zetas for curves defined over finite fields in the 
previous section. 
Main result  here is about a convergence region of such an Euler product. 
Key ingredients of our proof are a result of
(Harder-Narasimhan) Siegel, an ugly yet very precise formula 
for our local zeta functions,
Clifford Lemma for semi-stable vector bundles, and Weil's theorem 
on Riemann Hypothesis for Artin zeta functions.
\vskip 0.30cm
\centerline {\bf 2.1. Preparations}
\vskip 0.30cm
\noindent
{\bf 2.1.1. Invariants $\alpha,\beta$ and $\gamma$.} 
Let $C$ be a regular, reduced, irreducible
projective curve of genus $g$ defined over the finite field
${\Bbb F}_q$ with $q$  elements. As in I.1, we then get 
(the subset of ${\Bbb F}_q$-rational
points of) the associated moduli spaces
${\Cal M}_{E,r}(L)$ and ${\Cal M}_{C,r}(d)$. Recall that in I.1,
motivated by a work of Harder-Narasimhan [HN], we, following 
Desale-Ramanan [DR], defined the Harder-Narasimhan numbers 
$\beta_{C,r,{\Bbb F}_q}(L), \beta_{C,r,{\Bbb F}_q}(d)$, which
are very useful in the discussion of our zeta functions.  
Now  we introduce new invariants for $C$ by setting
$$\alpha_{C,r,{\Bbb F}_q}(d):=\sum_{V\in [V]\in {\Cal M}_{C,r}(d)({\Bbb
F}_q)}{{q^{h^0(C,V)}}\over {\#\text{Aut}(V)}},\qquad\qquad
\gamma_{C,r,{\Bbb F}_q}(d):=\sum_{V\in [V]\in {\Cal M}_{C,r}(d)({\Bbb
F}_q)}{{q^{h^0(C,V)}-1}\over {\#\text{Aut}(V)}},$$
and similarly define $\alpha_{C,r,{\Bbb F}_q}(L)$ and 
$\gamma_{C,r,{\Bbb F}_q}(L)$.
\vskip 0.30cm
\noindent
{\bf Lemma.} {\it With the same notation as above,

\noindent
(a) for $\alpha_{C,r,{\Bbb F}_q}(d)$,  $$\alpha_{C,r,{\Bbb F}_q}(d)=\cases
\beta_{C,r,{\Bbb F}_q}(d);&
d< 0;\\
\alpha_{C,r,{\Bbb F}_q}(r(2g-2)-d)\cdot q^{d-r(g-1)},& 0\leq d\leq r(2g-2);\\
\beta_{C,r,{\Bbb F}_q}(d)\cdot q^{d-r(g-1)},&
d>r(2g-2);\endcases$$

\noindent
(b) for $\beta_{C,r,{\Bbb F}_q}(d)$,
$$\beta_{C,r,{\Bbb F}_q}(\pm d+rn)=\beta_{C,r,{\Bbb F}_q}(d)
\qquad n\in {\Bbb Z};$$

\noindent
(c) for $\gamma_{C,r,{\Bbb F}_q}(d)$, $$\gamma_{C,r,{\Bbb
F}_q}(d)=\alpha_{C,r,{\Bbb F}_q}(d)-\beta_{C,r,{\Bbb F}_q}(d).$$}

\noindent
{\it Proof.}  (c) simply comes from the definition, while (b) is a direct 
consequence of Lemma 2 in 1.2.2 and
the fact that $\text{Aut}(V)\simeq \text{Aut}(V^\vee)$ 
for a vector bundle $V$. So it suffices to
prove (a).

When $d<0$, the relation is deduced from the fact that $h^0(C,V)=0$ if
$V$ is a semi-stable vector bundle with strictly negative degree; when 
$0\leq d\leq r(2g-2)$, the
result comes from the Riemann-Roch and Serre duality; finally when 
$d>r(2g-2)$, the result
is a direct consequence of the Riemann-Roch and the fact that  
$h^1(C,V)=0$ if $V$ is
a semi-stable vector bundle with degree strictly bigger  than $r(2g-2)$.
\vskip 0.30cm
We here reminder the reader that this Lemma and  Lemma 2 in 1.2.2 tell us 
that all $\alpha_{C,r,{\Bbb F}_q}(d),\beta_{C,r,{\Bbb F}_q}(d)$ and 
$\gamma_{C,r,{\Bbb F}_q}(d)$'s for all $d\in {\Bbb Z}$ may be
calculated from a finite subset of them, that is, from
$\alpha_{C,r,{\Bbb F}_q}(i),\beta_{C,r,{\Bbb F}_q}(j)$ 
with $i=0,\dots, r(g-1)$ and
$j=0,\dots,r-1$.
\vskip 0.30cm
\noindent
{\bf 2.1.2. Asymptotic Behaviors of $\alpha,\beta$ and $\gamma$}

For later use, we here discuss the asymptotic behavior of 
$\alpha_{C,r,{\Bbb F}_q}(d)$, 
$\beta_{C,r,{\Bbb F}_q}(d)$, and $\gamma_{C,r, {\Bbb F}_q}(0)$ 
when $q\to \infty$.
\vskip 0.30cm
\noindent
{\bf Proposition.} {\it With the same notation as above, when $q\to\infty$,

\noindent
(a) For all $d$, $$\beta_{C,r,{\Bbb
F}_q}(d)=O\Big(q^{r^2(g-1)}\Big);$$

\noindent
(b) $${{q^{(r-1)(g-1)}}\over{\gamma_{C,r,{\Bbb
F}_q}(0)}}=O\Big(1\Big).$$

\noindent
(c) For $0\leq d\leq r(g-1)$,
$${{\alpha_{C,r,{\Bbb F}_q}(d)}\over
{q^{d/2+r+r^2(g-1)}}}=O(1).$$}

\noindent
{\it Proof.}  Following Harder and Narasimhan [HN],  a result of Siegel 
on quadratic forms which is
equivalent to the fact that Tamagawa number of $\text{SL}_r$ is 1, 
may be understood via the
following relation on automorphism groups of rank $r$ vector bundles:
$$\sum_{V:r(V)=r,\text{det}(V)=L}{1\over {\#\text{Aut}(V)}}={{q^{(r^2-1)(g-1)}}
\over {q-1}}\cdot
\zeta_C(2)\dots\zeta_C(r).$$ Here $V$ runs over all rank $r$ vector 
bundles with determinant $L$ and
$\zeta_C(s)$ denotes the Artin zeta function of $C$.
Thus, $$0<\beta_{C,r,{\Bbb F}_q}(L)\leq {{q^{(r^2-1)(g-1)}}\over {q-1}}\cdot
\zeta_C(2)\dots\zeta_C(r).$$ This implies $$\beta_{C,r,{\Bbb
F}_q}(d)=\prod_{i=1}^{2g}(1-\omega_{C,1,{\Bbb F}_q}(i))\cdot 
\beta_{C,r,{\Bbb F}_q}(L)
\leq \prod_{i=1}^{2g}(1-\omega_{C,1,{\Bbb F}_q}(i))\cdot{{q^{(r^2-1)(g-1)}}
\over {q-1}}\cdot
\zeta_C(2)\dots\zeta_C(r).$$ Here two facts are used:

\noindent
(1) The number of ${\Bbb F}_q$-rational points of degree $d$ Jacobian 
$J^d(C)$ 
is equal to
$\prod_{i=1}^{2g}(1-\omega_{C,1,{\Bbb F}_q}(i))$; and

\noindent
(2) a result of Desale and Ramanan, which says that for any two 
$L,L'\in \text{Pic}^d(C)$,
$\beta_{C,r,{\Bbb F}_q}(L)=\beta_{C,r,{\Bbb F}_q}(L')$. 
(See e.g., [DR, Prop 1.7.(i)])

Thus by Weil's theorem on Riemann Hypothesis on Artin zeta functions  ([W1]),  
$$|\omega_{C,1,{\Bbb
F}_q}(i)|=O(q^{1/2}), \qquad i=0,\dots,2g.$$ This then completes the 
proof of (a).

To prove (b), first  note that (b) is equivalent to that, asymptotically, 
the lower bound of
$\gamma_{C,r,{\Bbb F}_q}(0)$ is at least $q^{(r-1)(g-1)}$. To show this, 
note that
$$\eqalign{\gamma_{C,r,{\Bbb F}_q}(0)&\geq
\sum_{V={\Cal O}_C\oplus L_2\oplus\dots\oplus L_r,
L_2,\dots,L_r\in \text{Pic}^0(C), \#\{{\Cal O}_C,L_2,\dots,L_r\}=r}
{{q^{h^0(C,V)}-1}\over {\#\text{
Aut}(V)}}\cr
 =&{1\over {(q-1)^{r-1}}}\sum_{V={\Cal O}_C\oplus L_2
\oplus\dots\oplus L_r,
L_2,\dots,L_r\in \text{Pic}^0(C), 
\#\{{\Cal O}_C,L_2,\dots,L_r\}=r} 1.\cr}$$ 
Now, by 
the above cited result of Weil again, as $q\to\infty$, 
$$\sum_{V={\Cal O}_C\oplus L_2\oplus\dots\oplus L_r,
L_2,\dots,L_r\in \text{Pic}^0(C), \#\{{\Cal O}_C,L_2,\dots,L_r\}=r} 
1=O(q^{g(r-1)}).$$
So we have (b) as well.

Just as (a), (c) is about to give an upper bound for 
$\alpha_{C,r,{\Bbb F}_q}(d)$ for $0\leq d\leq
r(2g-2)$. For this, we first recall the following

\noindent
{\bf Clifford Lemma.} (See e.g., [B-PBGN, Theorem
2.1]) {\it Let $V$ be a semi-stable bundle of rank
$r$ and degree
$d$ with $0\leq\mu(V)\leq 2g-2$. Then $$h^0(C,V)\leq r+{d\over 2}.$$}
 
Thus, $$\alpha_{C,r,{\Bbb F}_q}(d)\leq q^{{d\over 2}+r}\cdot
\beta_{C,r,{\Bbb F}_q}(d).$$
With this,  (c) is a direct consequence of (a).
\eject
\vskip 0.30cm
\noindent
{\bf 2.1.3. Ugly Formula}
\vskip 0.30cm
Recall that the rationality of $\zeta_{C,r,{\Bbb F}_q}(s)$ says that 
there exists a degree $2rg$
polynomial $P_{C,r,{\Bbb F}_q}(t)\in {\Bbb Q}[t]$ such that
$$Z_{C,r,{\Bbb F}_q}(t)={{P_{C,r,{\Bbb F}_q}(t)}
\over {(1-t^r)(1-q^rt^r)}}.$$
Thus we may set $$P_{C,r,{\Bbb F}_q}(t)=:
\sum_{i=0}^{2rg}a_{C,r,{\Bbb F}_q}(i)t^i.$$ On the other
hand, by the functional equation for $\xi_{C,r,{\Bbb F}_q}(t)(s)$, we have
$$P_{C,r,{\Bbb F}_q}(t)=
P_{C,r,{\Bbb F}_q}({1\over {qt}})\cdot q^{rg}\cdot t^{2rg}.$$ 
Hence, by comparing coefficients on both sides, we get the following
\vskip 0.30cm
\noindent
{\bf Lemma.} {\it With the same notation as above, for $i=0,1,\dots,rg-1$,
$$a_{C,r,{\Bbb F}_q}(2rg-i)=a_{C,r,{\Bbb F}_q}(i)\cdot q^{rg-i}.$$}

Now, to determine $P_{C,r,{\Bbb F}_q}(t)$ and hence 
$\zeta_{C,r,{\Bbb F}_q}(s)$ it suffices to find
$a_{C,r,{\Bbb F}_q}(i)$ for $i=0,1,\dots ,rg$.
\vskip 0.30cm
\noindent
{\bf Proposition.} ({\bf An Ugly Formula}) {\it With the same notation 
as above,
$$\eqalign{~&a_{C,r,{\Bbb F}_q}(i)\cr
=&\cases \alpha_{C,r,{\Bbb
F}_q}(d)-\beta_{C,r,{\Bbb F}_q}(d),& 0\leq i\leq r-1;\\
\alpha_{C,r,{\Bbb F}_q}(d)
-(q^r+1)\alpha_{C,r,{\Bbb F}_q}(d-r)+q^r\beta_{C,r,{\Bbb
F}_q}(d-r),& r\leq i\leq 2r-1;\\
\alpha_{C,r,{\Bbb F}_q}(d)
-(q^r+1)\alpha_{C,r,{\Bbb F}_q}(d-r)+q^r\alpha_{C,r,{\Bbb
F}_q}(d-2r),& 2r\leq i\leq r(g-1)-1;\\
-(q^r+1)\alpha_{C,r,{\Bbb F}_q}(r(g-2))
+q^r\alpha_{C,r,{\Bbb F}_q}(r(g-3))
+\alpha_{C,r,{\Bbb F}_q}(r(g-1)),& i=r(g-1);\\
\alpha_{C,r,{\Bbb F}_q}(d)
-(q^r+1)\alpha_{C,r,{\Bbb F}_q}(d-r)+\alpha_{C,r,{\Bbb
F}_q}(d-2r)q^r,& r(g-1)+1\leq i\leq rg-1;\\
2q^r\alpha_{C,r,{\Bbb F}_q}(r(g-2))
-(q^r+1)\alpha_{C,r,{\Bbb F}_q}(r(g-1)),& i=rg;\endcases
\cr}$$}

\noindent
{\it Proof.} By definition,
$$\eqalign{~&Z_{C,r,{\Bbb F}_q}(t)\cr
=&(\sum_{d=0}^{r(2g-2)}+\sum_{d=r(2g-2)+1}^\infty)\sum_{V\in [V]\in
{\Cal M}_{C,r}(d),d\geq 0} {{q^{h^0(C,V)}-1}\over{\#\text{Aut}(V)}}t^d\cr
=&\sum_{d=0}^{r(2g-2)}\sum_{V\in [V]\in {\Cal
M}_{C,r}(d)} {{q^{h^0(C,V)}-1}\over{\#\text{Aut}(V)}}t^d\cr
&\qquad+\sum_{i=1}^r\sum_{n=2g-2}^\infty\sum_{d=rn+i} \sum_{V\in [V]
\in {\Cal
M}_{C,r}(d)} {{q^{rn+i-r(g-1)}-1}\over{\#\text{Aut}(V)}}t^{rn+i}\cr
=&\sum_{d=0}^{r(2g-2)}\sum_{V\in [V]\in {\Cal
M}_{C,r}(d)({\Bbb F}_q)} {{q^{h^0(C,V)}-1}\over{\#\text{Aut}(V)}}t^d\cr
&\qquad+{{q^{r(1-g)}}\over{1-(qt)^r}}(qt)^{r(2g-2)}
 \sum_{i=1}^r \beta_{C,r,{\Bbb F}_q}(i) (qt)^i
-{1\over{1-t^r}}t^{r(2g-2)}\sum_{i=1}^r \beta_{C,r,{\Bbb F}_q}(i) t^i,\cr}$$
by a similar calculation as in the proof of Lemma 1.2.4.(b).
Now
$$\sum_{d=0}^{r(2g-2)}=\sum_{d=0,r(2g-2)}+\sum_{d=1,r(2g-2)-1}+\dots+
\sum_{d=r(g-1)-1,r(g-1)+1}+\sum_{d=r(g-1)}.$$
Thus, by Riemann-Roch,  Serre duality and Lemma 2.1.1, we conclude that
$$\eqalign{~&\sum_{d=0}^{r(2g-2)}\sum_{V\in [V]\in {\Cal
M}_{C,r}(d)} {{q^{h^0(C,V)}-1}\over{\#\text{Aut}(V)}}t^d\cr
=&\sum_{d=0}^{r(g-1)-1}\Big[\alpha_{C,r,{\Bbb F}_q}(d)
\Big(t^d+q^{r(g-1)-d}t^{r(2g-2)-d}\Big)-\beta_{C,r,{\Bbb
F}_q}(d)\Big(t^d+t^{r(2g-2)-d}\Big)\Big]\cr
&\qquad\qquad+\Big(\alpha_{C,r,{\Bbb
F}_q}(r(g-1))-\beta_{C,r,{\Bbb
F}_q}(r(g-1))\Big)\cdot t^{r(g-1)}.\cr}$$ With all this, together with 
Lemma 2 in 1.2.2
and the Lemma in 2.1.1, by a couple of pages routine calculation, we are 
lead to the 
ugly yet very precise formula in the proposition.
\vskip 0.30cm
\centerline {\bf 2.2. Global Non-Abelian Zeta Functions for Curves}
\vskip 0.30cm
\noindent
{\bf 2.2.1. Definition.} Let ${\Cal C}$ be a regular, reduced, irreducible 
projective curve of genus
$g$ defined over a number field $F$. Let $S_{\text{bad}}$ be the 
collection of all infinite places and
those finite places of $F$ at which ${\Cal C}$  does not have good 
reductions. As usual, a place $v$ of
$F$ is called good if $v\not\in S_{\text{bad}}$.
 
Thus, in particular, for any good place $v$ of $F$,   the $v$-reduction 
of ${\Cal C}$, denoted as
${\Cal C}_v$, gives a regular, reduced, irreducible projective curve 
defined over the residue 
field $F(v)$ of $F$ at $v$. Denote the cardinal number of $F(v)$ by 
$q_v$.  Then,  by the
construction of I.1, we obtain  associated rank $r$ 
non-abelian zeta function $\zeta_{{\Cal
C}_v,r,{\Bbb F}_{q_v}}(s)$. Moreover, from the rationality of 
$\zeta_{{\Cal
C}_v,r,{\Bbb F}_{q_v}}(s)$, there exists a degree $2rg$ polynomial 
$P_{{\Cal C}_v,r,{\Bbb F}_{q_v}}(t)\in {\Bbb Q}[t]$
such that
$$Z_{{\Cal C}_v,r,{\Bbb F}_{q_v}}(t)={{P_{{\Cal C}_v,r,{\Bbb F}_{q_v}}(t)}
\over {(1-t^r)(1-q^rt^r)}}.$$
Clearly, $$P_{{\Cal C}_v,r,{\Bbb F}_{q_v}}(0)
=\gamma_{{\Cal C}_v,r,{\Bbb F}_{q_v}}(0)\not=0.$$ Thus it
makes sense to introduce the polynomial 
$\tilde P_{{\Cal C}_v,r,{\Bbb F}_{q_v}}(t)$ with constant
term 1 by setting
$$\tilde P_{{\Cal C}_v,r,F(v)}(t):={{P_{{\Cal C}_v,r,F(v)}(t)}\over 
{P_{{\Cal C}_v,r,F(v)}(0)}}.$$
Now by definition, {\it the rank $r$ non-abelian 
zeta function $\zeta_{{\Cal C},r,F}(s)$ of} ${\Cal C}$ over $F$ is the 
following Euler product
$$\zeta_{{\Cal C},r,F}(s)
:=\prod_{v:\text{good}}{1\over{
\tilde P_{{\Cal C}_v,r,{\Bbb F}_{q_v}}(q_v^{-s})}},\hskip 2.0cm 
\text{Re}(s)>>0.$$

Clearly, when $r=1$, $\zeta_{{\Cal C},r,F}(s)$ coincides with the classical 
Hasse-Weil zeta function for $C$ over $F$ ([H]).
\vskip 0.30cm
\noindent
{\bf 2.2.2. Convergence.} At this earlier stage of the study of our 
non-abelian zeta functions, the
central problem is to justify the above definition. That is to say, 
to show  the above Euler product converges. In this direction, 
we have the following
\vskip 0.30cm
\noindent
{\bf Theorem.} {\it  Let ${\Cal C}$ be a regular, reduced, irreducible 
projective 
curve defined over
a number  field $F$. Then its associated rank r global non-abelian zeta 
function 
$\zeta_{{\Cal C},r,F}(s)$ converges  when $\text{Re}(s)\geq 1+g+(r^2-r)(g-1)$.}
\vskip 0.30cm
\noindent
{\it Proof.} Clearly, it suffices to show that for the reciprocal roots 
$\omega_{C,r,{\Bbb F}_q}(i),
i=1,\dots,2rg$, of $P_{C,r,{\Bbb F}_q}(t)$ associated to curves $C$ over
 finite fields ${\Bbb F}_q$,
$$|\omega_{C,r,{\Bbb F}_q}(i)|=O(q^{g+(r^2-1)(g-1)}).$$
Thus we are lead to estimate coefficients of $P_{C,r,{\Bbb F}_q}(t)$. 
Since we have the ugly yet very
precise formula for these coefficients, i.e., the Lemma and the 
Proposition in 2.1.3, 
it suffices to give upper
bounds for $\alpha_{C,r,{\Bbb F}_q}(i), \beta_{C,r,{\Bbb F}_q}(j)$ and a 
lower bound for $\gamma_{C,f,{\Bbb F}_q}(0)$,
 the constant term of $P_{C,r,{\Cal F}_q}(t)$. 
Thus, to complete the proof, we only need to cite the Proposition in 2.1.2. 
\vskip 0.30cm
\noindent
{\bf Question.} {\it For any regular, reduced, irreducible projective curve 
${\Cal C}$ of genus $g$ defined over a number field $F$, whether  its 
associated rank $r$ global non-abelian zeta function
$\zeta_{{\Cal C},r,F}(s)$  admits  meromorphic continuation to 
the whole complex $s$-plane.}
\vskip 0.30cm
Recall that even  when $r=1$, i.e., for the classical Hasse-Weil 
zeta functions, 
this is still quite open. 
\vskip 0.30cm
\noindent
{\bf 2.2.3. Working Hypothesis.} Like in the theory of abelian zeta functions, 
we want to use our non-abelian zeta functions  to study  non-abelian 
arithmetic aspect of curves. Motivated by the classical analytic class 
number formula for Dedekind zeta functions and its counterpart BSD 
conjecture for Hasse-Weil zeta functions of elliptic curves, we 
expect that our non-abelian zeta function can be used to understand 
the Weil-Petersson volumes of
moduli space of stable bundles as well as the associated Tamagawa measures.

As such, local factors for \lq bad' places are needed. Our suggestion is 
as follows: 
for $\Gamma$-factors,  we take those coming from the functional equation for
 $\zeta_F(rs)\cdot\zeta_F(r(s-1))$, where $\zeta_F(s)$ denotes
the standard Dedekind zeta function of $F$; while for finite bad places, 
we  first use the semi-stable reduction for curves to find a semi-stable model 
for ${\Cal C}$, then use Seshadri's moduli spaces of parabolic bundles 
to construct polynomials for singular fibers, which
usually have degree lower than $2rg$. With all this done, we 
then can introduce the so-called completed rank
$r$ non-abelian zeta function  for ${\Cal C}$ over
$F$, or better,  the  completed rank $r$ non-abelian zeta 
function $\xi_{X,r,{\Cal O}_F}(s)$ for a
semi-stable model $X\to \text{Spec}({\Cal O}_F)$ of ${\Cal C}$. 
Here ${\Cal O}_F$ denotes the ring of
integers of $F$. (If necessary, we take a finite extension of $F$.)

\noindent
{\bf Question.} {\it Whether the meromorphic extension of
$\xi_{X,r,{\Cal O}_F}(s)$, if exists,  satisfies the functional equation
$$\xi_{X,r,{\Cal O}_F}(s)=\pm\, \xi_{X,r,{\Cal O}_F}(1+{1\over r}-s).$$}
{\it Remark.} From our study for non-abelian zeta functions of elliptic curves 
[We3], we obtain the
following  \lq absolute Euler product' for rank 2 zeta 
functions of elliptic curves
$$\eqalign{\zeta_2(s)=&\prod_{p>2; \text{prime}}{1\over
{1+(p-1)p^{-s}+(2p-4)p^{-2s}+(p^2-p)p^{-3s}+p^2p^{-4s}}}\cr
=&\prod_{p>2; \text{prime}}{1\over {A_p(s)+B_p(s)p^{-2s}}},
\qquad \text{Re}(s)>2\cr}$$
with $$A_p(s)=1+(p-1)p^{-s}+(p-2)p^{-2s},\qquad B_p(s)
=(p-2)+(p^2-p)p^{-s}+p^2p^{-2s}.$$
Set $t:=q^{-s}$ and $a_p(t):=A_p(s), b_p(t):=B_p(s)$. Then in 
${\Bbb Z}[t]$, we have the factorization
$$a_p(t)=(1+(p-2)t)(1+t),\qquad
b_p(t)=((p-2)+pt)(1+pt)$$ and
$$a_p({1\over {pt}})={1\over {p^2t^2}}\cdot b_p(t).$$ As pointed to me by
Kohnen, $$1+(p-1)p^{-s}+(2p-4)p^{-2s}+(p^2-p)p^{-3s}+p^2p^{-4s}$$
 is quite similar to Andrianov's genus two spinor 
$L$-function. (See e.g. [We1].)
\vskip 0.45cm
\centerline {\li I.3. Non-Abelian Zeta
Functions and Infinitesimal Structures of Brill-Noether Loci}
\vskip 0.30cm
In this section, we study the infinitesimal structures of the so-called 
non-abelian Brill-Noether
loci for rank two semi-stable vector bundles over genus two curves. As an 
application, we calculate
the corresponding rank two non-abelian zeta functions. During this process, 
we
see clearly how Weierstrass points, intrinsic arithmetic invariants of 
curves [We1], contribute to
our zeta functions among others. 

We in this section assume that the characteristic of the base field is 
strictly bigger than 2 for
simplicity.
\vskip 0.30cm
\centerline {\bf 3.1. Invariants $\beta_{C,2, {\Bbb F}_q}(d)$.} 

Let $C$ be a genus two
regular reduced irreducible
projective curve defined over ${\Bbb F}_q$. Here we want to calculate 
Harder-Narasimhan numbers
$\beta_{C,2,{\Bbb F}_q}(d)$ for all $d$. Note that from the Lemma  in 2.1.1, 
$$\beta_{C,2,{\Bbb F}_q}(d)
=\beta_{C,2,{\Bbb F}_q}(d+2n).$$ So it suffices to calculate 
$\beta_{C,2,{\Bbb F}_q}(d)$ when $d=0,1$.
For this, we cite the following result of Desale and Ramanan:

\noindent
{\bf Proposition.} ([DR]) {\it With the same notation as above, for 
$L\in \text{Pic}^d(C), d=0,1$,
$$\beta_{C,2,{\Bbb F}_q}(L)={{q^3}\over
{q-1}}\cdot\zeta_C(2)-q\prod_{i=1}^4(1-\omega_i)\cdot 
\sum_{d_1+d_2=d, d_1>d_2}{{\beta_{C,1,{\Bbb
F}_q}(d_1)\beta_{C,1,{\Bbb F}_q}(d_2)}\over {q^{d_1-d_2}}}.$$ Here 
$\zeta_C(s)$ denotes the Artin
zeta function for $C$ and $\omega_1,\dots,\omega_4$ are the roots of 
the associated
$Z$-function $Z_C(s)$, i.e., $\omega_{C,1,{\Bbb F}_q}(i), i=0,\dots, 
4=2\times 2$ 
in our notation.}

Thus, in particular, $\beta_{C,2,{\Bbb F}_q}(L)$ is independent of $L$. 

\noindent
{\bf Lemma.} {\it With the same notation as above, for $d=0,1$
$$\beta_{C,2,{\Bbb F}_q}(d)={{q^3}\over {q-1}}\cdot\zeta_C(2)\cdot 
\prod_{i=1}^4(1-\omega_i)
-{{q^{d+1}}\over {(q-1)^2(q^2-1)}}\cdot \prod_{i=1}^4(1-\omega_i)^4.$$}

\noindent
{\it Proof.} This comes from the following two facts:

\noindent
(1) for all $d$, $$\beta_{C,1,{\Bbb F}_q}(d)=
{{\prod_{i=1}^4(1-\omega_i)}\over{q-1}};$$ 

\noindent
(2) the number of ${\Bbb F}_q$-rational points of $\text{Pic}^d(C)$ is equal to
$\prod_{i=1}^4(1-\omega_i)$.
\vskip 0.30cm
\centerline {\bf 3.2. Invariants $\alpha\ \&\  \gamma$: Easy Parts}
\vskip 0.30cm
\noindent
{\bf 3.2.1. Infinitesimal structures: a taste.} Here we want to calculate 
$\alpha_{C,2,{\Bbb
F}_q}(0)$. By the Lemma in 3.1, it suffices to give 
$\gamma_{C,2,{\Bbb F}_q}(0)$. So we are lead to
study
$\gamma_{C,2,{\Bbb F}_q}(L)$ which is
supported over the Brill-Noether locus $$W_{C,2}^0(L):=\{[V]\in 
{\Cal M}_{C,2}(L):h^0(C,\text{
Gr}(V))\geq 1\}.$$ (In general, as in  [B-PGN], 
we define the Brill-Noether locus by
$$W_{C,2}^k(L):=\{[V]\in {\Cal M}_{C,2}(L):h^0(C,\text{Gr}(V))\geq k+1\}.)$$

Note that no degree zero stable bundle admits non-trivial
global sections, so
$W_{C,2}^0(L):=\{[{\Cal O}_C\oplus L]\}$ consists of only one single point.

\noindent
(a) If $L={\Cal O}_C$, then $W_{C,2}^0({\Cal O}_C)=W_{C,2}^1({\Cal O}_C).$
 Moreover, infinitesimally,
$V={\Cal O}_C\oplus {\Cal O}_C$ or $V$ corresponds to all non-trivial 
extensions 
$$0\to {\Cal O}_C\to V\to {\Cal O}_C\to 0$$ which are parametrized by 
${\Bbb P}\text{Ext}^1({\Cal
O}_C,{\Cal O}_C)\simeq {\Bbb P}^1$. Thus, by definition,
$$\gamma_{C,2,{\Bbb F}_q}({\Cal O}_C)={{q^2-1}\over {(q^2-1)(q^2-q)}}+
(q+1)\cdot {{q-1}\over {q(q-1)}}
={q\over {q-1}}.$$

\noindent
(b) If $L\not ={\Cal O}_C$, then, 
infinitesimally,
$V={\Cal O}_C\oplus L$ or $V$ corresponds to the single non-trivial extension 
$$0\to {\Cal O}_C\to V\to L\to 0.$$  Thus, by definition,
$$\gamma_{C,2,{\Bbb F}_q}(L)={{q-1}\over {(q-1)^2}}+{{q-1}\over {q-1}}
={q\over {q-1}}.$$

Thus we have the following

\noindent
{\bf Lemma.} {\it With the same notation as above, for all 
$L\in \text{Pic}^0(C)$,
$$\gamma_{C,2,{\Bbb F}_q}(L)={q\over {q-1}}.$$  In particular, 
$$\gamma_{C,2,{\Bbb F}_q}(0)={q\over {q-1}}\cdot \prod_{i=1}^4(1-\omega_i).$$}

\noindent
{\bf 3.2.2. Invariants $\alpha_{C,2{\Bbb F}_q}(1)$.} As before, 
it suffices to calculate
$\gamma_{C,2,{\Bbb F}_q}(L)$ for all $L\in \text{Pic}^1(C)$. Note that 
in this case, all bundles are
stable, so $\text{Aut}(V)\simeq {\Bbb F}_q^*$ and
$$W_{C,2}^0(L)\simeq\{V:\text{stable}, r(V)=2,\text{det}(V)=L, h^0(C,V)\geq
1\}.$$ Moreover, by [B-PGN, Prop. 3.1], $$W_{C,2}^0(L)=\{V:\text{stable}, 
r(V)=2,\text{
det}(V)=L,h^0(C,V)=1\}$$ and any $V\in W_{C,2}^0(L)$ admits a non-trivial 
extension
$$0\to {\Cal O}_C\to V\to L\to 0.$$ On the other hand, any non-trivial 
extension
$$0\to {\Cal O}_C\to V\to L\to 0$$ gives rise to a stable bundle. 
So in fact
$$W_{C,2}^0(L)\simeq {\Bbb P}\text{Ext}^1(L,{\Cal O}_C)\simeq {\Bbb P}^1.$$
Thus we have the following

\noindent
{\bf Lemma.} {\it With the same notation as above, for $L\in \text{Pic}^1(C)$,
$$W_{C,2}^0(L)\simeq {\Bbb P}^1,\qquad \text{and}\qquad 
\gamma_{C,2,{\Bbb F}_q}(L)=q+1.$$
In particular, $$\gamma_{C,2{\Bbb F}_q}(1)=(q+1)\cdot 
\prod_{i=1}^4(1-\omega_i).$$}
\vskip 0.30cm
\centerline {\bf 3.3. Infinitesimal Structures of Non-Abelian Brill-Noether 
Loci}
\vskip 0.30cm 
We next calculate $\gamma_{C,r,{\Bbb F}_q}(2)$. In general, 
the level $r(g-1)$, 
which in our present case corresponds to 2, is
the most complicated one. So  the discussion here is rather involved.

Let us start with the structures of the non-abelian Brill-Noether loci 
$W_{C,2}^0(L)$ and $W_{C,2}^1(L)$ for $L\in \text{Pic}^2(C)$.
For this,  recall the structure  map 
$\pi:C\times C/S_2\to \text{Pic}^2(C)$ defined by $[(x,y)]\mapsto 
{\Cal O}_C(x+y)$. 
Here $S_2$
denotes the symmetric group of two symbols which acts naturally on 
$C\times C$ via $(x,y)\mapsto
(y,x)$. One checks that $\pi$ is a one point blowing-up centered at 
the canonical line
bundle $K_C$ of $C$. For later use, denote by $\Delta$ the image of 
the diagonal of $C\times C$ in
$\text{Pic}^2(C)$.

Next, we want to understand the structure of sublocus 
$W_{C,2}^0(L)^{\text{ss}}$ of $W_{C,2}^0(L)$
consisting of non-stable but semi-stable vector bundles.
By definition, for any $V\in [V]\in  W_{C,2}^0(L)^{\text{ss}}$, 
$\text{Gr}(V)={\Cal O}_C(P)\oplus L(-P)$
for a suitable (${\Bbb F}_q$-rational) point $P\in C$. 
Thus accordingly,

\noindent
(a) if $L\not=K_C$, then $W_{C,2}^0(L)^{\text{ss}}$ is parametrized by 
(${\Bbb F}_q$-rational points
of) $C$, due to the fact that now $h^0(C,L)=1$. Write also 
$L={\Cal O}_C(A+B)$ with two points $A,B$
of $C$, which are unique from the above discussion on the map $\pi$, we 
then conclude that
$$W_{C,2}^1(L)=\{[{\Cal O}_C(A)\oplus {\Cal O}_C(B)]\}.$$ 

\noindent
(b) if $L=K_C$, then for any $P$, $K_C={\Cal O}_C(P+\iota(P))$ where 
$\iota: C\to C$ denotes the
canonical involution on $C$. So $$W_{C,2}^0(L)^{\text{ss}}
=\{[{\Cal O}_C(P)\oplus {\Cal O}_C(\iota(P))]:P\in C\}.$$ Therefore
 $W_{C,2}^0(L)^{\text{ss}}$ is parametrized by ${\Bbb
P}^1$. Moreover, $$W_{C,2}^1(K_C)=W_{C,2}^0(L)^{\text{ss}}=\{[{\Cal O}_C(P)
\oplus {\Cal
O}_C(\iota(P))]:P\in C\}.$$

On the other hand, it is easy to check that every non-trivial extension 
$$0\to {\Cal O}_V\to W\to
L\to 0$$ gives rise to a semi-stable vector bundle $W$, and if $W$ is not 
stable, then there exists
a point $Q\in C$ such that $W$ may also be given by the non-trivial 
extension $$0\to {\Cal O}_C(Q)\to
W\to L(-Q)\to 0.$$ Note also that the kernel of the natural map 
$H^1(C,\text{Hom}(L,{\Cal
O}_C))\to H^1(C,\text{Hom}(L(-Q),{\Cal O}_C))$ is one dimensional. 
So among all
non-trivial extensions $0\to {\Cal O}_C\to V\to L\to 0$, which are 
parametrized by ${\Bbb P}\text{
Ext}^1(L,{\Cal O}_C)\simeq {\Bbb P}^2$, the non-stable (yet semi-stable) 
vector bundles are
parametrized by (${\Bbb F}_q$-rational points of) $C$ when $L\not=K_C$ by (a), 
or ${\Bbb P}^1$ when
$L=K_C$ by (b). (See [NR, Lemma 3.1]) In this way, we have proved the
following  result on non-abelian Brill-Noether loci for moduli space of
${\Cal M}_{C,2}(L)$ with $L$ a degree 2 line bundle on a genus two curve, 
which is not covered by
[B-PGN]:
\vskip 0.30cm
\noindent
{\bf Lemma.} {\it With the same notation as above, 
$W_{C,2}^0(L)\simeq {\Bbb P}\text{Ext}^1(L,{\Cal
O}_C)\simeq {\Bbb P}^2,$ in which the locus $W_{C,2}^0(L)^{\text{ss}}$ 
of semi-stable but not stable
bundles is parametrized by
$C$ or ${\Bbb P}^1$  according to  $L\not= K_C$ or $L=K_C$. More precisely,

\noindent
(a) if $L={\Cal O}_C(A+B)\not=K_C$ with $A,B$ two points of $C$, then 
$W_{C,2}^0(L)^{\text{ss}}$, as a
birational image of $C$ under the complete linear system $K_C(A+B)$, 
is a degree 4 plane curve 
with a single node located at $W_{C,2}^1(L)=\{[{\Cal O}_C(A)\oplus 
{\Cal O}_C(B)]\}$;

\noindent
(b) If $L=K_C$, as a degree 2 regular plane curve,
$$W_{C,2}^1(K_C)=W_{C,2}^0(L)^{\text{ss}}=\{[{\Cal O}_C(P)\oplus 
{\Cal O}_C(\iota(P))]:P\in C\}\simeq
{\Bbb P}^1.$$} 

Next, we study the infinitesimal structures of non-abelian Brill-Noether 
loci. Set 
$$W_{C,2}^0(L)^{s}:=W_{C,2}^0(L)\backslash W_{C,2}^0(L)^{\text{ss}}.$$ 
Then the infinitesimal structure of $W_{C,2}^0(L)$ at points 
$[V]\in W_{C,2}^0(L)^{s}$ is simple:
each $[V]$ consists a single stable rank two vector bundle with 
$\text{det}(V)=L, h^0(C,V)=1$ and
$\text{Aut}(V)\simeq {\Bbb F}_q^*$.

Now we consider $W_{C,2}^0(L)^{\text{ss}}$. 

\noindent
(a)  $L\not=K_C$. Then there exist two points $A,B$ of $C$ such that
 $L={\Cal O}_C(A+B)$. 
Thus, for any $V\in [{\Cal O}_C(P)\oplus {\Cal O}_C(A+B-P)]\not\in 
W_{C,2}^1(L)$,
$V$ is given by an extension $0\to {\Cal O}_C(P)\to V\to 
{\Cal O}_C(A+B-P)\to 0$ due to the fact that
for the non-trivial extension $0\to {\Cal O}_C(A+B-P)\to W\to 
{\Cal O}_C(P)\to 0$, $h^0(C,W)=0$.
Thus,  each class $[{\Cal O}_C(P)\oplus {\Cal O}_C(A+B-P)]\not\in W_{C,2}^1(L)$
consists of exactly two vector bundles, i.e., 
$V_1={\Cal O}_C(P)\oplus {\Cal O}_C(A+B-P)$ and $V_2$ given by 
the non-trivial extension 
$0\to {\Cal O}_C(P)\to V\to {\Cal O}_C(A+B-P)\to 0$. Clearly, 
$h^0(C,V_1)=h^0(C,V_2)=1$ and $\#\text{
Aut}(V_1)=(q-1)^2,  \#\text{Aut}(V_2)=q-1;$
\vskip 0.30cm
To study $W_{C,2}^1(L)=\{[{\Cal O}_C(A)\oplus {\Cal O}_C(B)]\}$, we 
divide it into  two subcases.

\noindent
(i)  $A\not=B$. Then there are exactly three vector bundles $V_0,\, V_1$ 
and $V_2$
in the class $[{\Cal O}_C(A)\oplus {\Cal
O}_C(B)]$. They are $V_0={\Cal O}_C(A)\oplus {\Cal O}_C(B)$, 
$V_1$ given by the non-trivial extension $0\to {\Cal O}_C(A)\to V_1\to 
{\Cal O}_C(B)\to 0$ and $V_2$ 
given by the non-trivial extension $0\to {\Cal O}_C(B)\to V_2\to 
{\Cal O}_C(A)\to 0$. Clearly,
$h^0(C,V_1)=2, h^0(C,V_1)=h^0(C,V_2)=1$ and $\#\text{Aut}(V_0)=(q-1)^2,\ 
\#\text{Aut}(V_1)=\#\text{Aut}(V_2)=q-1$;

Thus in particular, $$\gamma_{C,2,{\Bbb F}_q}(L)=(q^2+q+1-(N_1-1))
\cdot {{q-1}\over {q-1}}
+(N_1-2)\Big({{q-1}\over {(q-1)^2}}+{{q-1}\over {q-1}}\Big)+\Big({{q^2-1}
\over {(q-1)^2}}+
{{q-1}\over {q-1}}+{{q-1}\over {q-1}}\Big).$$ Here $N_1=q+1-
(\omega_1+\dots+\omega_4)$ denotes the
number of ${\Bbb F}_q$-rational points of $C$.

\noindent
(ii) $A=B$. Then the infinitesimal structure at $[{\Cal O}_C(A)\oplus 
{\Cal O}_C(A)]$ is as follows: an independent point corresponding to 
$V_0={\Cal O}_C(A)\oplus {\Cal O}_C(A)$ and a projective line 
parametrizing all non-trivial extension $0\to {\Cal O}_C(A)\to V\to
{\Cal O}_C(A)\to 0$. Clearly,
$h^0(C,V_0)=2, h^0(C,V)=1$ and $\#\text{Aut}(V_0)=(q^2-1)(q^2-q),\ 
\#\text{Aut}(V)=q(q-1);$

Thus in particular, $$\gamma_{C,2,{\Bbb F}_q}(L)=(q^2+q+1-(N_1-1))\cdot 
{{q-1}\over {q-1}}
+(N_1-2)\Big({{q-1}\over {(q-1)^2}}+{{q-1}\over {q-1}}\Big)+\Big({{q^2-1}
\over{(q^2-1)(q^2-q)}}+(q+1){{q-1}\over {q(q-1)}}\Big).$$

\noindent
(b) $L=K_C$. Then $K_C={\Cal O}_C(P+\iota(P))$ for all points $P$. 
Therefore for all $[V]\in 
W_{C,2}^1(L)=W_{C,2}^0(L)^{\text{ss}}$, $[V]=[{\Cal O}_C(P)\oplus 
{\Cal O}_C(\iota(P))]$. Accordingly,
two subcases:

\noindent
(i) $P\not=\iota P$. Then there are exactly three vector bundles 
$V_0,\, V_1$ and $V_2$
in the class $[{\Cal O}_C(P)\oplus
{\Cal O}_C(\iota(P)]$. They are $V_0={\Cal O}_C(P)\oplus 
{\Cal O}_C(\iota(P))$, 
$V_1$ given by the non-trivial extension $0\to {\Cal O}_C(P)\to V_1\to 
{\Cal O}_C(\iota(P))\to 0$
and $V_2$  given by the non-trivial extension $0\to 
{\Cal O}_C(\iota(P))\to V_2\to {\Cal O}_C(P)\to
0$. Clearly,
$h^0(C,V_0)=2, h^0(C,V_1)=h^0(C,V_2)=1$ and $\#\text{Aut}(V_0)=(q-1)^2,\ 
\#\text{Aut}(V_1)=\#\text{Aut}(V_2)=q-1$;

\noindent
(ii) $P=\iota(P)$ a Weierstrass point, all of which are six. Then the 
infinitesimal structure at
$[{\Cal O}_C(P)\oplus {\Cal O}_C(P)]$ is as follows: an independent point 
corresponding to $V_0={\Cal
O}_C(P)\oplus {\Cal O}_C(P)$ and a projective line parametrizing all 
non-trivial extension $0\to
{\Cal O}_C(P)\to V\to {\Cal O}_C(P)\to 0$. Clearly,
$h^0(C,V_0)=2, h^0(C,V)=1$ and $\#\text{Aut}(V_0)=(q^2-1)(q^2-q),\ 
\#\text{Aut}(V)=q(q-1).$

Thus, in particular, 
$$\eqalign{~&\gamma_{C,2,{\Bbb F}_q}(K_C)=(q^2+q+1-(q+1))\cdot {{q-1}\over
{q-1}}\cr
&\qquad+(q+1-6)\Big({{q^2-1}\over {(q-1)^2}}+
{{q-1}\over {q-1}}+{{q-1}\over {q-1}}\Big)+6\Big({{q^2-1}\over 
{(q^2-1)(q^2-q)}}+(q+1){{q-1}\over
{q(q-1)}}\Big).\cr}$$

All in all, we have completed the proof of the following

\noindent
{\bf Proposition.} {\it With the same notation as above,

\noindent
(a) For $L\not=K_C$, 

\noindent
(i) if $L\not\in \Delta$, $\gamma_{C,2,{\Bbb F}_q}(L)={{q^3+2q-3+N_1}
\over {q-1}}$;

\noindent
(ii) if $L\in \Delta$, $\gamma_{C,2,{\Bbb F}_q}(L)={{q^3-2+N_1}\over {q-1}}$;

\noindent
(b) if $L=K_C$, $\gamma_{C,2,{\Bbb F}_q}(L)={{q^3+2q^2-10q+5}\over{q-1}}$.

\noindent
In particular, 
$$\gamma_{C,2,{\Bbb F}_q}(2)=\Big(\prod_{i=1}^4(1-\omega_i)-(q+1)\Big)\cdot
{{q^3+2q-3+N_1}\over {q-1}}+q\cdot {{q^3-2+N_1}\over {q-1}}
+{{q^3+2q^2-10q+5}\over{q-1}}.$$}

In this way, by using the ugly formula in  2.1.3, we can finally write down 
the rank two
non-abelian zeta functions for genus two curves, where degree 8 polynomials 
are involved. We leave this to the reader.

\vskip 1.0cm
\centerline {\li Chapter II. Non-Abelian L-Functions}
\vskip 0.30cm
While we may introduce general non-abelian $L$-functions by 
using more general test functions as what Tate did in his Thesis where abelian 
version is discussed, in this paper, we decide to take a different approach 
using Eisenstein series. (We reminder the reader that for the abelian picture, 
Eisenstein series are not available.) Moreover, as for the integration domain,
we use a much more general type of moduli spaces.
\vskip 0.45cm
\centerline 
{\li II.1. Epstein Zeta Functions and Non-Abelian Zeta Functions}
\vskip 0.30cm
To motivate what follows, we begin this chapter with a discussion 
on non-ablian zeta functions for number fields.

For simplicity, assume that the  number field involved is the field of 
rationals. A lattice $\Lambda$ over ${\Bbb Q}$ is   semi-stable,
by definition, if for any sublattice $\Lambda_1$ of $\Lambda$,
$$\big(\text{Vol}\,\Lambda_1\big)^{\text{rank}\,\Lambda}\geq 
\big(\text{Vol}\,\Lambda\big)^{\text{rank}\,\Lambda_1}.$$
Denote the moduli space of rank $r$ semistable lattices over ${\Bbb Q}$ by 
${\Cal  M}_{{\Bbb Q},r}$, then the lattice version of {\it
rank $r$ non-abelian zeta function} $\xi_{{\Bbb Q},r}(s)$ of ${\Bbb Q}$ is
defined to be  
$$\xi_{{\Bbb Q},r}(s):=
\int_{{\Cal  M}_{{\Bbb Q},r}}\left( e^{h^0({\Bbb Q},\Lambda)}-1\right)
\cdot \big(e^{-s}\big)^{\text{deg}(\Lambda)}\, d\mu(\Lambda),\qquad 
\text{Re}(s)>1,$$
 where $h^0({\Bbb Q},\Lambda):=\log\left(\sum_{x\in \Lambda}
\exp\big(-\pi|x|^2\big)\right)$ and 
$\text{deg}(\Lambda):=
-\log\big(\text{Vol}({\Bbb R}^{\text{rank}(\Lambda)}/\Lambda)\big)$
denotes the Arakelov degree of $\Lambda$. 
Moreover, note that
the newly defined $h^0$ has a natural company $h^1$ and that similarly as
cohomology for bundles over curves, $h^i$ satisfy the Serre duality and
Riemann-Roch (for details, see [We2]). In particular, as shown in [We2], 
(see also the calculation below for an alternative proof,)

\noindent
(i) $\xi_{{\Bbb Q},1}(s)$ coincides with the (completed)
Riemann-zeta function;

\noindent
(ii) $\xi_{{\Bbb Q},r}(s)$ can be meromorphically extended to the whole 
complex plane;

\noindent
(iii) $\xi_{{\Bbb Q},r}(s)$ satisfies
the functional equation $$\xi_{{\Bbb Q},r}(s)=\xi_{{\Bbb Q},r}(1-s);$$ 

\noindent
(iv) $\xi_{{\Bbb Q},r}(s)$ has only two singularities, simple poles, 
at $s=0,1$, with the same residues
$\text{Vol}\left({\Cal  M}_{{\Bbb Q},r}[1]\right)$, the Tamagawa type 
volume of 
the space of rank $r$ semi-stable lattice of volume 1. 

Denote by ${\Cal  M}_{{\Bbb Q},r}[T]$ the moduli space of rank $r$ semi-stable 
lattices of volume $T$. We have a 
trivial decomposition $${\Cal  M}_{{\Bbb Q},r}=
\cup_{T>0}{\Cal  M}_{{\Bbb Q},r}[T].$$ Moreover, there is a natural morphism
$${\Cal  M}_{{\Bbb Q},r}[T]\to {\Cal  M}_{{\Bbb Q},r}[1],\qquad \Lambda\mapsto 
T^{1\over r}\cdot\Lambda.$$

With this, for $\text{Re}(s)>1$,
$$\align \xi_{{\Bbb Q},r}(s)=&
\int_{\cup_{T>0}{\Cal  M}_{{\Bbb Q},r}[T]}\left( e^{h^0({\Bbb Q},\Lambda)}-1
\right)\cdot \big(e^{-s}\big)^{\text{deg}(\Lambda)}\, d\mu(\Lambda)\\
=&\int_0^\infty T^s{{dT}\over T}
\int_{{\Cal  M}_{{\Bbb Q},r}[1]}\left( 
e^{h^0({\Bbb Q},T^{1\over r}\cdot \Lambda)}-1\right)
\cdot d\mu_1(\Lambda),\endalign$$ where $d\mu_1$ denotes the induced 
Tamagawa measure 
on ${\Cal  M}_{{\Bbb Q},r}[1]$.

Thus note that $$h^0({\Bbb Q},T^{1\over r}\cdot \Lambda)=
\log\left(\sum_{x\in \Lambda}
\exp\big(-\pi|x|^2\cdot T^{2\over r}\big)\right),$$
and for $B\not=0$,
$$\int_0^\infty e^{-AT^B}T^s{{dT}\over T}={1\over B}\cdot A^{-{s\over B}}
\cdot\Gamma({s\over B}),$$ we have
$$\xi_{{\Bbb Q},r}(s)={r\over 2}\cdot\pi^{-{r\over 2}\, s}
\Gamma({r\over 2}\, s)\cdot\int_{{\Cal  M}_{{\Bbb Q},r}[1]}
\left(\sum_{x\in\Lambda\backslash\{0\}}|x|^{-rs}\right)\cdot d\mu_1(\Lambda).$$
Set now the completed Epstein zeta function, a special kind of Eisenstein 
series, associated to the rank $r$ lattice 
$\Lambda$ over ${\Bbb Q}$ by
$$\hat E(\Lambda;s):=\pi^{-s}\Gamma(s)\cdot \sum_{x\in \Lambda\backslash \{0\}}
|x|^{-2s},$$ then we have the following
\vskip 0.30cm
\noindent
{\bf Proposition.} (Eisenstein series and Non-Abelian 
Zeta Functions) {\it With the same notation as above,
$$\xi_{{\Bbb Q},r}(s)={r\over 2}\int_{{\Cal  M}_{{\Bbb Q},r}[1]}\hat 
E(\Lambda,{r\over 2}s)\,d\mu_1(\Lambda).$$}

\noindent
{\it Remark.} Such a non-abelian zeta is indeed very beautiful: not only its 
construction is so elegent, its structure is also very rational -- 
Recently, Lagarias and Suzuki ([LS]) show that the rank two zeta 
$\xi_{{\Bbb Q},2}(s)$ for the field of rationals satisfies the Riemann 
Hypothesis, i.e., the zeros are all on the line 
$\Re(s)=\frac{1}{2}$. 
\vskip 0.20cm
\centerline
{\li II.2.   Canonical Polygons and Geometric Truncation}

We start with Weil's adelic interpretation of locally free sheaves on curves.
Fix a smooth geometrically connected projective curve $X$ over a finite field 
${\Bbb F}_q$. 
Denote its function field by $F$ and identify the places of $F$ with the 
closed points of 
$X$ which we denote by $|X|$. For each place $x$ of $F$, set $F_x$ the 
$x$-completion 
of $F$
 with ${\Cal O}_x$ the ring of integers,
$\pi_x$ a local parameter, and $\kappa(x)$ the residue field. Denote by 
$x:F_x^*\to {\Bbb Z}$ 
the normalized valuation of $F_x$ such that $x(\pi_x)=1$. Denote also by 
${\Bbb A}$ the ring
 of adeles and ${\Cal O}_{\Bbb A}$ the ring of integers.

If $E$ is  a locally free ${\Cal O}_F$-sheaf of rank $r$ over $X$, denote by 
$E_F$ 
the fiber 
of $E$ at the generic point $\text {Spec}(F)$ of $X$ ($E_F$ is an $F$-vector 
space of dimension $r$), 
and for each $v\in |X|$, set $E_{{\Cal O}_v}:=
H^0(\text {Spec}\,{\Cal O}_{F_v},E)$ a free ${\Cal O}_v$-module of rank 
$r$. In particular, we have a canonical isomorphism:
$$\text{can}_v:F_v\otimes_{{\Cal O}_v}E_{{\Cal O}_v}\simeq F_v\otimes_FE_F.$$
Thus, in particular,  with respect to a basis 
$\alpha_F:F^r\simeq E_F$ of its generic fiber and a basis 
$\alpha_{{\Cal O}_v}:{\Cal O}_v^r\simeq E_{{\Cal O}_v}$ for any 
$v\in |X|$, the elements $g_v:=(F_v\otimes_F\alpha_F)^{-1}\circ\text{can}_v
\circ (F_v\otimes_{{\Cal O}_v}\alpha_{{\Cal O}_v})\in GL_r(F_v)$ for all 
$v\in |X|$ define an element $g_{\Bbb A}:=(g_v)_{v\in |X|}$ of 
$GL_r({\Bbb A})$, since for almost all $v$ we have $g_v\in 
GL_r({\Cal O}_v)$. As a result, we obtain  a bijection from 
the set of isomorphism classes of triples $(E;\alpha_F;
(\alpha_{{\Cal O}_v})_{v\in |X|})$ as above onto $GL_r({\Bbb A})$. 
Moreover, if $r\in GL_r(F), k\in GL_r({\Cal O}_F)$ and if this bijection
 maps the triple $(E;\alpha_F;(\alpha_{{\Cal O}_v})_{v\in |X|})$ onto 
$g_{\Bbb A}$, the same map maps the triple $(E;\alpha_F\circ r^{-1};
(\alpha_{{\Cal O}_v}\circ k_v)_{v\in |X|})$ onto $rg_{\Bbb A}k$. 
Therefore the above bijection induces a bijection between the set of 
isomorphism classes of locally free ${\Cal O}_F$-sheaves of rank $r$ over $X$
 and the double coset space
$GL_r(F)\backslash GL_r({\Bbb A})/GL_r({\Cal O}_F)$.

More generally, let $r=r_1+\cdots+r_s$ be a partition $I=(r_1,\cdots,r_s)$ of 
$r$ and let $P_I$ be the corresponding standard parabolic subgroup of $GL_r$. 
Then we have a natural bijection from the set of isomorphism classes of 
triple $(E_*;\alpha_{*,F}:(\alpha_{*,{\Cal O_v}})_{v\in |X|})$ onto 
$P_I({\Bbb A})$, where $E_*:=\big((0)=E_0\subset E_1\subset\cdots\subset
 E_s\big)$ is a filtration of locally free sheaves of rank $(r_1,r_1+r_2,\cdots, 
r_1+r_2+\cdots+r_s=r)$ over $X$, (i.e, each $E_j$ is a vector sheaf of rank 
$r_1+r_2+\cdots+r_j$ over $X$ and each quotient $E_j/E_{j-1}$ is torsion 
free,) which is equipped with an isomorphism of filtrations of $F$-vector 
spaces
$$\alpha_{*,F}:\big((0)=F_0\subset F^{r_1}\subset\cdots\subset F^{r_1+r_2+
\cdots+r_s=r}\big)\simeq (E_*)_F,$$ and with an isomorphism of filtrations 
of free ${\Cal O}_v$-modules
$$\alpha_{*,{\Cal O}_v}:\big((0)\subset {\Cal O}_v^{r_1}\subset\cdots\subset 
{\Cal O}_v^{r_1+r_2+\cdots+r_s=r}\big)\simeq (E_*)_{{\Cal O}_v},$$ for every 
$v\in |X|$. Moreover this bijection induces a bijection between the set of 
isomorphism classes of the filtrations of locally free sheaves of rank 
$(r_1,r_1+r_2,\cdots, r_1+r_2+\cdots+r_s=r)$ over $X$ and the double coset 
space $P_I(F)\backslash P_I({\Bbb A})/P_I({\Cal O}_{\Bbb A})$. The natural 
embedding $P_I({\Bbb A})\hookrightarrow P_I({\Bbb A})$ (resp. the 
canonical projection $P_I({\Bbb A})\to M_I({\Bbb A})\to GL_{r_j}
({\Bbb A})$ for $j=1,\cdots,s$, where $M_I$ denotes the standard Levi 
of $P_I$) admits the modular interpretation
$$(E_*;\alpha_{*,F}:(\alpha_{*,{\Cal O_v}})_{v\in |X|})\mapsto 
(E_s;\alpha_{s,F}:(\alpha_{s,{\Cal O_v}})_{v\in |X|})$$ (resp. 
$$(E_*;\alpha_{*,F}:(\alpha_{*,{\Cal O_v}})_{v\in |X|})\mapsto 
(\text{gr}_j(E_*);\text{gr}_j(\alpha_{*,F}), \text{gr}_j
(\alpha_{*,{\Cal O_v}})_{v\in |X|}),$$ where $\text{gr}_j(E_*):=
E_j/E_{j-1}$,
$\text{gr}_j(\alpha_{*,F}):F^{r_j}\simeq \text{gr}_j(E_*)_F$ and 
$\text{gr}_j(\alpha_{*,{\Cal O_v}}):{\Cal O_v}^{r_j}\simeq \text{gr}_j
(E_*)_{\Cal O_v}$, $v\in |X|$ are induced by $\alpha_{*,F}$ and 
$\alpha_{*,{\Cal O_v}}$ respectively.)

Denote by $E_g$ the rank $r$ locally free sheaf on $X$ associated to
 $g\in GL_r({\Bbb A})$. Then, 
 $$\text{deg}(E_g)=-\log\big(N(\text{det}g)\big)$$ 
with $N:GL_1({\Bbb A}_F)={\Bbb I}_F\to {\Bbb Q}_{>0}$ the standard norm 
map of the idelic group of $F$.

With this, for $g\in GL_r({\Bbb A})$ and a parabolic subgroup 
$Q$ of $GL_r$, denote by $E_*^{g;Q}$ the filtration of the locally free sheaf 
$E_{g}$ induced by the parabolic subgroup $Q$. 

Now following Lafforgue [Laf], introduce an associated polygon 
$p_Q^g:[0,r]\to {\Bbb Q}$ by the following 3 conditions:

\noindent
(i) $p_Q^g(0)=p_Q^g(r)=0$;

\noindent
(ii) $p_Q^g$ is affine on the interval 
$[\text{rank}E_{i-1}^{g;Q},\text{rank}E_{i}^{g;Q}]$; and 

\noindent
(iii) for all indices $i$, 
$$p_Q^g(\text{rank}E_{i}^{g;Q})=\text{deg}(E_{i}^{g;Q},\rho_{i}^{g;Q})-
{{\text{rank}E_{i}^{g;Q}}\over r}\cdot \text{deg}(E_{g},\rho_{g}).$$
Then by Prop. 1 in I.1.1, i.e, the existence and uniqueness of 
Harder-Narasimhan filtration,
 there is a unique convex polygon $\bar p^g$ 
which bounds all $p_Q^g$ from above for all parabolic subgroups $Q$ for 
$GL_r$. 
Moreover there exists a parabolic subgroup $\bar Q^g$ such that 
$p_{{\bar Q}^g}^g=\bar p^g$.  In particular, as a direct consequence, we 
obtain the following
well-known 
\vskip 0.30cm
\noindent
{\bf Lemma.} (See e.g. [Laf]) {\it For any fixed polygon $p:[0,r]\to {\Bbb Q}$
 and any $d\in {\Bbb Z}$, the subset
$$\{g\in GL_r(F)\backslash GL_r({\Bbb A}):\text{deg}\,g=d,\bar p^g\leq p\}$$ 
is compact.}
\vskip 0.30cm 
Similarly yet more generally, for a fixed parabolic subgroup $P$ of $GL_r$ 
and $g\in GL_r({\Bbb A})$, there is a unique maximal element $\bar p_P^g$ 
among all $p_Q^g$, where $Q$ runs over all parabolic subgroups of $GL_r$ 
which are contained in $P$.  And we have
\vskip 0.30cm
\noindent
{\bf Lemma}$'$. (See e.g. [Laf]) 
{\it For any fixed polygon $p:[0,r]\to {\Bbb Q}$, $d\in {\Bbb Z}$ 
and  any standard parabolic subgroup $P$ of $GL_r$, the subset
$$\{g\in GL_r(F)\backslash GL_r({\Bbb A}):\text{deg}\,g=d,\bar p^g_P\leq p, 
p^g_P\geq -p\}$$ is compact.}
\vskip 0.30cm
Moreover, let $p,q:[0,r]\to {\Bbb R}$ be two polygons and $P$ a standard 
parabolic 
subgroup of $GL_r$. Then as in [Laf], we say $q>_Pp$ if for any $1\leq i\leq 
|P|$,
$$q(\text{rank}E_i^P)>p(\text{rank}E_i^P)$$ where $(r_1,\cdots,r_{|P|})$ 
denotes the partition of $r$ corresponding to $P$. As usual denote by 
${\bold 1}$ the characteristic function of the variable $g\in GL_r({\Bbb A})$. 
For example,
$${\bold 1}(\bar p^g\leq p)(g)=\cases 1,&\text {if}\ p^g\leq p\\
0&\text{otherwise}.\endcases$$
Then we have the following result of Lafforgue.
\vskip 0.30cm
\noindent
{\bf Proposition.} ([Laf, Prop. V.1.c]) {\it For any convex polygon 
$p:[0,r]\to {\Bbb R}$, as a 
function of $g\in GL_r({\Bbb A})$,
$${\bold 1}(\bar p^g\leq p)=\sum_{P\supset P_0}(-1)^{|P|-1}
\sum_{\delta\in P(F)\backslash GL_r(F)}{\bold 1}(p_P^{\delta g}>_Pp).$$ 
Here $P$ runs over all standard parabolic subgpoups of $GL_r$.}
\vskip 0.60cm
\centerline
{\li II.3. Non-Abelian L-Functions}
\vskip 0.30cm
In this section, we introduce non-abelian $L$-functions for function fields 
and study their basic properties.
\vskip 0.30cm
\centerline {\bf 3.1. Choice of Moduli Spaces}

For the function field $F$ with genus $g_X$, and for a fixed 
$r\in {\Bbb Z}_{>0}$, we take the moduli space to be 
$${\Cal M}^{\leq p}_{F,r}:=\{g\in GL_r(F)Z_{GL_r({\Bbb A})}\backslash 
GL_r({\Bbb A}):\bar p^g\leq p\}$$ 
for a fixed convex polygon $p:[0,r]\to{\Bbb R}$. Also 
we denote by $d\mu$ the induced Tamagawa measures on 
${\Cal M}_{F,r}^{\leq p}$. 

More generally, for any standard parabolic subgroup $P$ of $GL_r$, we 
introduce the moduli spaces $${\Cal M}_{F,r}^{P;\leq p}
:=\{g\in P(F)Z_{GL_r({\Bbb A})}\backslash GL_r({\Bbb A}):
\bar p_P^g\leq p, \bar p_P^g\geq -p\}.$$ By the discussion in II.2, these 
moduli spaces 
${\Cal M}_{F,r}^{P;\leq p}$ are all compact, a key
 property which plays a central role in our definition of non-abelian 
$L$-functions below.
\vfill
\eject
\vskip 0.30cm
\centerline {\bf 3.2. Choice of Eisenstein Series: 
First Approach to Non-Abelian $L$-Function}

To facilitate our ensuing discussion, we start with  some preparations. 
For details,
please consult [MW], which is heavily used in this subsection. (The 
experienced
reader may skip this subsection, except for possible later reference 
about notations.)

Fix a connected reduction group $G$ defined over $F$, denote by $Z_G$ its 
center. Fix a minimal parabolic subgroup $P_0$ of $G$. Then $P_0=M_0U_0$, 
where as usual we fix once and for all the Levi $M_0$ and  the unipotent
 radical $U_0$. A parabolic subgroup $P$ of $G$ is called standard if 
$P\supset P_0$. For such groups write $P=MU$ with $M_0\subset M$ the 
standard Levi and $U$ the unipotent radical. Denote by $\text {Rat}(M)$ 
the group of rational characters of $M$, i.e, the morphism $M\to {\Bbb G}_m$ 
where ${\Bbb G}_m$ denotes the multiplicative group. Set 
$$\frak a_M^*:=\text {Rat}(M)\otimes_{\Bbb Z}{\Bbb C},\qquad \frak a_M
:=\text{Hom}_{\Bbb Z}(\text {Rat}(M),{\Bbb C}),$$ and $$\text {Re}\,\frak 
a_M^*:=\text {Rat}(M)\otimes_{\Bbb Z}{\Bbb R},\qquad \text{Re}\,\frak a_M
:=\text{Hom}_{\Bbb Z}(\text {Rat}(M),{\Bbb R}).$$ For any $\chi\in 
\text {Rat}(M)$, we obtain a (real) character $|\chi|:M({\Bbb A})\to 
{\Bbb R}^*$ defined by $m=(m_v)\mapsto m^{|\chi|}:=\prod_{v\in S}
|m_v|_v^{\chi_v}$ with $|\cdot|_v$ the $v$-absolute values. Set then 
$M({\Bbb A})^1:=\cap_{\chi\in \text {Rat}(M)}\text{Ker}|\chi|$, which is 
a normal subgroup of $M({\Bbb A})$. Set $X_M$ to be the group of complex 
characters which are trivial on $M({\Bbb A})^1$. Denote by 
$H_M:=\log_M:M({\Bbb A})\to \frak a_M$ the map such that 
$\forall\chi\in \text {Rat}(M)\subset \frak a_M^*,\langle\chi,
\log_M(m)\rangle:=\log(m^{|\chi|})$. Clearly, $$M({\Bbb A})^1
=\text{Ker}(\log_M);\qquad \log_M(M({\Bbb A})/M({\Bbb A})^1)\simeq 
\text{Re}\,\frak a_M.$$ Hence in particular there is a natural isomorphism 
$\kappa:\frak a_M^*\simeq X_M.$
Set $$\text{Re}\,X_M:=\kappa (\text{Re}\,\frak a_M^*),\qquad \text{Im}\,X_M:=
\kappa (i\cdot \text{Re}\,\frak a_M^*).$$ Moreover define our working space 
$X_M^G$ to be the subgroup of $X_M$ consisting of complex characters of 
$M({\Bbb A})/M({\Bbb A})^1$ which are trivial on $Z_{G({\Bbb A})}$.

Fix a maximal compact subgroup ${\Bbb K}$ such that for all standard parabolic
subgroups $P=MU$ as above, $P({\Bbb A})\cap{\Bbb K}=M({\Bbb A})
\cap{\Bbb K}\cdot
U({\Bbb A})\cap{\Bbb K}.$ Hence we get the Langlands decomposition 
$G({\Bbb A})=M({\Bbb A})\cdot U({\Bbb A})\cdot {\Bbb K}$. Denote by 
$m_P:G({\Bbb A})\to M({\Bbb A})/M({\Bbb A})^1$ the map $g=m\cdot n\cdot
 k\mapsto M({\Bbb A})^1\cdot m$ where $g\in G({\Bbb A}), m\in M({\Bbb A}), 
n\in U({\Bbb A})$ and 
$k\in {\Bbb K}$.
 
Fix Haar measures on $M_0({\Bbb A}), U_0({\Bbb A}), {\Bbb K}$ respectively 
such that 

\noindent
(1) the induced measure on $M(F)$ is the counting measure and the volume of
 the induced measure on $M(F)\backslash M({\Bbb A})^1$ is 1. (Recall that 
it is  a fundamental fact that $M(F)\backslash M({\Bbb A})^1$ is compact.)

\noindent
(2) the induced measure on $U_0(F)$ is the counting measure and the volume 
of $U(F)\backslash U_0({\Bbb A})$ is 1. (Recall that being unipotent radical, 
$U(F)\backslash U_0({\Bbb A})$ is compact.)

\noindent
(3) the volume of ${\Bbb K}$ is 1.

Such measures then also induce Haar measures via $\log_M$ to $\frak a_{M_0}, 
\frak a_{M_0}^*$, etc. Furthermore, if we denote by $\rho_0$ a half of the 
sum of the positive roots of  the maximal split torus $T_0$ of the central 
$Z_{M_0}$
of $M_0$, then $$f\mapsto \int_{M_0({\Bbb A})\cdot U_0({\Bbb A})\cdot 
{\Bbb K}}f(mnk)\,dk\,dn\,m^{-2\rho_0}dm$$ defined for continuous functions 
with compact supports on $G({\Bbb A})$  defines a Haar measure $dg$ on 
$G({\Bbb A})$. This in turn gives measures on $M({\Bbb A}), U({\Bbb A})$ 
and hence on $\frak a_{M}, \frak a_{M}^*$, $P({\Bbb A})$, etc, for all 
parabolic subgroups $P$. In particular, one checks that the following 
compactibility condition holds
$$\int_{M_0({\Bbb A})\cdot U_0({\Bbb A})\cdot {\Bbb K}}f(mnk)\,dk\,dn\,
m^{-2\rho_0}dm=\int_{M({\Bbb A})\cdot U({\Bbb A})\cdot {\Bbb K}}f(mnk)\,dk\,
dn\,m^{-2\rho_P}dm$$ for all continuous functions $f$ with compact supports 
on $G({\Bbb A})$, where $\rho_P$ denotes a half of the sum of the positive 
roots of  the maximal split torus $T_P$ of the central $Z_{M}$ of $M$. 
For later use, 
denote also by $\Delta_P$ the set of positive roots determined by $(P,T_P)$, 
$\Delta_0=\Delta_{P_0}$ and $W$ the associated Weyl group.

Fix an isomorphism $T_0\simeq {\Bbb G}_m^R$ and a place $v_0$ of $F$ and a 
uniformizer
$\pi_{v_0}$ at $v_0$. The group $\pi_{v_0}^{\Bbb Z}$ generated by $\pi_{v_0}$ 
can be identified
with a subgroup of ${\Bbb A}^*$ and hence $(\pi_{v_0}^{\Bbb Z})^R$ 
with a subgroup of $T_0({\Bbb A})$. Thus there exists a $W$-invariant
subgroup of $Z_{M_0({\Bbb A})}$ 
which is isomorphic to a subgroup of finite index of $(\pi_{v_0}^{\Bbb Z})^R$. 
Fix such a group once and for all and denote it by $A_{M_0({\Bbb A})}$. 

 More generally, for a standard parabolic subgroup 
$P=MU$, set $A_{M({\Bbb A})}:=
A_{M_0({\Bbb A})}\cap Z_{M({\Bbb A})}$ where as used above $Z_*$ denotes 
the center of the group $*$. Then $A_{M({\Bbb A})}\backslash M({\Bbb A})/ 
M({\Bbb A})^1$ is finite. For later use, set also 
$A_{M({\Bbb A})}^G:=\{a\in A_{M({\Bbb A})}:\log_Ga=0\}.$ 
Then $A_{M({\Bbb A})}$ contains $A_{G({\Bbb A})}\oplus
A_{M({\Bbb A})}^G$ as a subgroup of finite index.

Note that ${\Bbb K}$, $M(F)\backslash M({\Bbb A})^1$ and $U(F)\backslash 
U({\Bbb A})$ are all compact, thus with the Langlands decomposition 
$G({\Bbb A})=U({\Bbb A})M({\Bbb A}){\Bbb K}$ in mind, the reduction theory 
for $G(F)\backslash G({\Bbb A})$ or more generally 
$P(F)\backslash G({\Bbb A})$ is reduced to that for 
$A_{M({\Bbb A})}$. As such for $t_0\in M_0({\Bbb A})$, set
$$A_{M_0({\Bbb A})}(t_0):=\{a\in A_{M_0({\Bbb A})}:a^\alpha>t_0^\alpha\ 
\forall\alpha\in\Delta_0\}.$$ Then, for a fixed compact subset 
$\omega\subset P_0({\Bbb A})$, we have the corresponding Siegel set 
$$S(\omega;t_0):=\{p\cdot a\cdot k:p\in \omega, a\in A_{M_0({\Bbb A})}
(t_0),k\in {\Bbb K}\}.$$ In particular, for big enough $\omega$ and small 
enough $t_0$, i.e, $t_0^\alpha$ is very close to 0 for all 
$\alpha\in\Delta_0$, the classical reduction theory may be restated as 
$G({\Bbb A})=G(F)\cdot S(\omega;t_0)$. More generally set 
$$A_{M_0({\Bbb A})}^P(t_0):=\{a\in A_{M_0({\Bbb A})}:a^\alpha>t_0^\alpha,\ 
\forall\alpha\in\Delta_0^P\},$$ and $$S^P(\omega;t_0):=\{p\cdot a\cdot k:
p\in \omega, a\in A_{M_0({\Bbb A})}^P(t_0),k\in {\Bbb K}\}.$$  
Then similarly as above
for big enough $\omega$ and small enough $t_0$, $G({\Bbb A})=P(F)\cdot 
S^P(\omega;t_0)$. (Here 
$\Delta_0^P$ denotes the set of positive roots for $(P_0\cap M,T_0)$.)

Fix an embedding $i_G:G\hookrightarrow SL_n$ sending $g$ to $(g_{ij})$. 
Introducing a hight function on $G({\Bbb A})$ by setting 
$\|g\|:=\prod_{v\in S}\text{sup}\{|g_{ij}|_v:\forall i,j\}$. It is well-known 
that up to $O(1)$, hight functions are unique. This implies that the following 
growth conditions do not depend on the height function we choose.

A function $f:G({\Bbb A})\to {\Bbb C}$ is said to have moderate growth if 
 there exist 
$c,r\in {\Bbb R}$ such that $|f(g)|\leq c\cdot \|g\|^r$ for all
 $g\in G({\Bbb A})$. Similarly, for a standard parabolic subgroup $P=MU$, 
a function $f:U({\Bbb A})M(F)\backslash G({\Bbb A})\to {\Bbb C}$ is said 
to have moderate growth if there exist $c,r\in {\Bbb R},\lambda\in 
\text{Re}X_{M_0}$ such that for any $a\in A_{M({\Bbb A})},k\in {\Bbb K}, 
m\in M({\Bbb A})^1\cap S^P(\omega;t_0)$,
$$|f(amk)|\leq c\cdot \|a\|^r\cdot m_{P_0}(m)^\lambda.$$

Now fix a place $v_0$ of $F$, denote by $G({\Bbb A})_{v_0}$ the 
inverse image of $G(F_{v_0})$ in $G({\Bbb A})$.
Denote by $\frak z$ the Bernstein centre of $G({\Bbb A})_{v_0}$. 
The $\frak z$ acts naturally on the locally constant functions on 
$G({\Bbb A})$.

By definition,  a function 
$\phi:U({\Bbb A})M(F)\backslash G({\Bbb A})\to {\Bbb C}$ is called automorphic 
if

\noindent
(i) $\phi$ has moderate growth;

\noindent
(ii) $\phi$ is locally constant;

\noindent
(iii) $\phi$ is ${\Bbb K}$-finite, i.e, the ${\Bbb C}$-span of all 
$\phi(k_1\cdot *\cdot k_2)$ parametrized by $(k_1,k_2)\in {\Bbb K}\times 
{\Bbb K}$ is finite dimensional; 

\noindent
(iv) $\phi$ is $\frak z$-finite, i.e, the ${\Bbb C}$-span of all 
$\delta(X)\phi$ parametrized by all $X\in \frak z$ is finite dimensional.

For such a function $\phi$, set $\phi_k:M(F)\backslash M({\Bbb A})\to 
{\Bbb C}$ to be the function defined by $m\mapsto m^{-\rho_P}\phi(mk)$ for all
 $k\in {\Bbb K}$. 
Set $A(U({\Bbb A})M(F)\backslash G({\Bbb A}))$ be the space of automorphic 
forms on $U({\Bbb A})M(F)\backslash G({\Bbb A})$.

For a measurable locally $L^1$-function $f:U(F)\backslash G({\Bbb A})\to 
{\Bbb C}$, define its constant term along with the standard parabolic subgroup 
$P=UM$ to be the function $f_P:U({\Bbb A})\backslash G({\Bbb A})\to {\Bbb C}$ 
given by
$g\to\int_{U(F)\backslash G({\Bbb A})}f(ng)dn.$ 
Then an automorphic form $\phi\in A(U({\Bbb A})M(F)\backslash G({\Bbb A}))$ is 
called a cusp form if for any standard parabolic subgroup $P'$ properly 
contained in $P$, $\phi_{P'}\equiv 0$. Denote by
 $A_0(U({\Bbb A})M(F)\backslash G({\Bbb A}))$  the space of cusp 
forms on $U({\Bbb A})M(F)\backslash G({\Bbb A})$. One checks easily that

\noindent
(i) all cusp forms are rapidly decreasing; and hence

\noindent
(ii) there is a natural pairing
 $$\langle\cdot,\cdot\rangle:A_0(U({\Bbb A})M(F)\backslash G({\Bbb A}))\times 
A(U({\Bbb A})M(F)\backslash G({\Bbb A}))\to {\Bbb C}$$ defined by
$\langle \psi,\phi\rangle:=\int_{Z_{M({\Bbb A})}U({\Bbb A})M(F)\backslash 
G({\Bbb A})}\psi(g)\bar\phi(g)\,dg.$

Moreover, for a (complex) character $\xi:Z_{M({\Bbb A})}\to {\Bbb C}^*$ of 
$Z_{M({\Bbb A})}$ set 
$$\eqalign{A&(U({\Bbb A})M(F)\backslash G({\Bbb A}))_\xi\cr
:=&\{\phi\in A(U({\Bbb A})M(F)
\backslash G({\Bbb A})):\phi(zg)=z^{\rho_P}\cdot\xi(z)\cdot\phi(g),\forall 
z\in Z_{M({\Bbb A})}, g\in G({\Bbb A})\}\cr}$$ and 
$$A_0(U({\Bbb A})M(F)\backslash G({\Bbb A}))_\xi:=A_0(U({\Bbb A})M(F)
\backslash G({\Bbb A}))\cap A(U({\Bbb A})M(F)\backslash G({\Bbb A}))_\xi.$$

Set now $$A(U({\Bbb A})M(F)\backslash G({\Bbb A}))_Z:=\sum_{\xi\in 
\text{Hom}(Z_{M({\Bbb A})},{\Bbb C}^*)}A(U({\Bbb A})M(F)\backslash 
G({\Bbb A}))_\xi$$ and
$$A_0(U({\Bbb A})M(F)\backslash G({\Bbb A}))_Z:=\sum_{\xi\in 
\text{Hom}(Z_{M({\Bbb A})},{\Bbb C}^*)}A_0(U({\Bbb A})M(F)\backslash 
G({\Bbb A}))_\xi.$$ It is well-known that the natural morphism $${\Bbb C}
[\text{Re}\frak a_M]\otimes A(U({\Bbb A})M(F)\backslash G({\Bbb A}))_Z
\to A(U({\Bbb A})M(F)\backslash G({\Bbb A}))$$ defined by 
$(Q,\phi)\mapsto \big(g\mapsto Q(\log_M(m_P(g))\big)\cdot \phi(g)$ is an 
isomorphism, using the special structure of $A_{M({\Bbb A})}$-finite 
functions and the Fourier analysis over the compact space $A_{M({\Bbb A})}
\backslash Z_{M({\Bbb A})}$. Consequently, we also obtain a natural isomorphism
$${\Bbb C}[\text{Re}\frak a_M]\otimes A_0(U({\Bbb A})M(F)\backslash 
G({\Bbb A}))_Z\to A_0(U({\Bbb A})M(F)\backslash G({\Bbb A}))_\xi.$$

Set also $\Pi_0(M({\Bbb A}))_\xi$ be isomorphism classes of irreducible 
representations of $M({\Bbb A})$ occurring in the space $A_0(M(F)\backslash 
M({\Bbb A}))_\xi$, and
$$\Pi_0(M({\Bbb A})):=\cup_{\xi\in \text{Hom}(Z_{M({\Bbb A})},{\Bbb C}^*)}
\Pi_0(M({\Bbb A}))_\xi.$$ For any 
$\pi\in \Pi_0(M({\Bbb A}))_\xi$, set $A_0(M(F)\backslash M({\Bbb A}))_\pi$ to 
be the isotypic component of type $\pi$ of $A_0(M(F)\backslash 
M({\Bbb A}))_\xi$, i.e, the set of cusp forms of $M({\Bbb A})$ generating a 
semi-simple isotypic $M({\Bbb A})$-module of type $\pi$.
Set $$\eqalign{A_0&(U({\Bbb A})M(F)\backslash G({\Bbb A}))_\pi\cr
:=&\{\phi\in 
A_0(U({\Bbb A})M(F)\backslash G({\Bbb A})):\phi_k\in A_0(M(F)\backslash 
M({\Bbb A}))_\pi,\forall k\in {\Bbb K}\}.\cr}$$ Clearly
$$A_0(U({\Bbb A})M(F)\backslash G({\Bbb A}))_\xi=\oplus_{\pi\in 
\Pi_0(M({\Bbb A}))_\xi} A_0(U({\Bbb A})M(F)\backslash G({\Bbb A}))_\pi.$$

More generally, let $V\subset A(M(F)\backslash M({\Bbb A}))$ be an irreducible 
$M({\Bbb A})$-module with $\pi_0$ the induced representation of 
$M({\Bbb A})$. Then we call $\pi_0$ an automorphic 
representation of $M({\Bbb A})$. Denote by $A(M(F)\backslash 
M({\Bbb A}))_{\pi_0}$  the isotypic 
subquotient module of type $\pi_0$ of $A(M(F)\backslash M({\Bbb A}))$. Then
$$V\otimes \text{Hom}_{M({\Bbb A})}(V,A(M(F)
\backslash M({\Bbb A})))\simeq A(M(F)\backslash M({\Bbb A}))_{\pi_0}.$$ 
Set $$\eqalign{A&(U({\Bbb A})M(F)\backslash G({\Bbb A}))_{\pi_0}\cr
:=&\{\phi\in 
A(U({\Bbb A})M(F)\backslash G({\Bbb A})):\phi_k\in A(M(F)\backslash 
M({\Bbb A}))_{\pi_0},\forall k\in {\Bbb K}\}.\cr}$$
Moreover if 
$A(M(F)\backslash M({\Bbb A}))_{\pi_0}\subset
A_0(M(F)\backslash M({\Bbb A}))$, we call $\pi_0$ a cuspidal representation.

Two automorphic representations $\pi$ and $\pi_0$ of $M({\Bbb A})$ are said 
to be equivalent if there exists $\lambda\in X_M^G$ such that $\pi\simeq 
\pi_0\otimes\lambda$. This, in practice, means that $A(M(F)\backslash 
M({\Bbb A}))_\pi
=\lambda\cdot A(M(F)\backslash M({\Bbb A}))_{\pi_0}.$ That is, for any 
$\phi_\pi\in A(M(F)\backslash M({\Bbb A}))_\pi$ there exists a 
$\phi_{\pi_0}\in A(M(F)\backslash M({\Bbb A}))_{\pi_0}$ such that 
$\phi_\pi(m)=m^\lambda\cdot \phi_{\pi_0}(m)$. Consequently, 
$$A(U({\Bbb A})M(F)\backslash G({\Bbb A}))_\pi=(\lambda\circ m_P)\cdot 
A(U({\Bbb A})M(F)\backslash G({\Bbb A}))_{\pi_0}.$$ Denote by 
$\frak P:=[\pi_0]$ the equivalence class of $\pi_0$. Then $\frak P$ is an 
$X_M^G$-principal homogeneous space, hence admits a natural complex structure. 
Usually we call $(M,\frak P)$ a cuspidal datum of $G$ if $\pi_0$ is cuspidal. 
Also for $\pi\in\frak P$ set $\text{Re}\,\pi:=\text{Re}\,\chi_\pi=|\chi_\pi|
\in 
\text{Re}X_M$, where $\chi_\pi$ is the central character of $\pi$, and 
$\text{Im}\,\pi:=\pi\otimes(-\text{Re}\,\pi)$.

Now fix an irreducible automorphic representation $\pi$ of $M({\Bbb A})$ and
$\phi\in A(U({\Bbb A})M(F)\backslash G({\Bbb A}))_\pi$, define the 
associated Eisenstein series $E(\phi,\pi):G(F)\backslash G({\Bbb A})\to 
{\Bbb C}$ by $$E(\phi,\pi)(g):=\sum_{\delta\in P(F)\backslash G(F)}\phi
(\delta g).$$ Then  there is an open cone ${\Cal C}\subset 
\text{Re}X_M^G$ such that if $\text{Re}\,\pi\in {\Cal C}$, $E(\lambda\cdot 
\phi,
\pi\otimes\lambda)(g)$ converges uniformly for $g$ in a compact subset of 
$G({\Bbb A})$ and $\lambda$ in an open neighborhood of 0 in $X_M^G$. For 
example, if $\frak P=[\pi]$ is cuspidal, we may even take ${\Cal C}$ to be 
the cone $\{\lambda\in \text{Re}X_M^G:\langle\lambda-\rho_P,
\alpha^\vee\rangle>0,\forall\alpha\in \Delta_P^G\}$.
As a direct consequence,  $E(\phi,\pi)\in A(G(F)\backslash G({\Bbb A}))$.
 That is, Eisenstein series $E(\phi,\pi)$ are automorphic forms.

As noticed above, being an automorphic form, $E(\phi,\pi)$ is of 
moderate growth.
 However, in general it is not integrable over $Z_{G({\Bbb A})}G(F)\backslash 
G({\Bbb A})$. To remedy this, classically, as initiated in the so-called 
Rankin-Selberg method, analytic truncation is used: From Fourier analysis, 
we understand 
that the problematic terms are the so-called constant terms, which are of 
moderate growth, so by cutting off them, the reminding one is 
rapidly increasing and hence integrable.

In general, it is very difficult to make such an analytic truncation 
intrinsically related with arithmetic properties of number fields. 
(See however,  [Z] and  [Ar1,2].) 
On the other hand, Eisenstein series themselves are quite intrinsic 
arithmetical invariants.
Thus it is natural for us on one hand to keep Eisenstein series 
unchanged while on the other 
to find new moduli spaces, which themselves intrinsically parametrize 
certain modular objects, and over which Eisenstein series are integrable.

This is exactly what we are doing now. As said, we are going to view Eisenstein
 series as something globally defined, and use a geometric truncation 
for the space $G(F)\backslash G({\Bbb A})$ so that the integrations of the 
Eisenstein series over the newly obtained moduli spaces give us naturally 
non-abelian $L$ functions for function fields.

As such, let us now come back to the group $G=GL_r$, then as in 3.1, we 
obtain the moduli space ${\Cal M}_{F,r}^{\leq p}$ 
and hence a well-defined integration $$L_{F,r}^{\leq p}(\phi,\pi):=
\int_{{\Cal M}_{F,r}^{\leq p}}E(\phi,\pi)(g)\,dg,
\qquad \text {Re}\,\pi\in {\Cal C}.$$
\vskip 0.10cm
\centerline
{\bf 3.3. New Non-Abelian $L$-Functions}
\vskip 0.30cm
However,  in such a general form, we do not know whether the latest  
defined integration has any nice properties, such as meromorphic continuation 
and functional equations, etc... It is to remedy this that 
we make a further selection about automorphic forms.
  
Fix then a convex polygon $p:[0,r]\to {\Bbb R}$ as in II.2 so as to 
obtain the moduli space ${\Cal M}_{F,r}^{\leq p}$. 
Set $G=GL_r$, fix the minimal parabolic subgroup $P_0$ corresponding to the 
partition $(1,\cdots,1)$ with $M_0$ consisting of diagonal matrices. Fix a 
standard parabolic subgroup $P_I=U_IM_I$ corresponding to the partition 
$I=(r_1,\cdots,r_{|P|})$ of $r$ with $M_I$ the standard Levi and $U_I$ the 
unipotent radical.
 
Then for a fixed irreducible automorphic representation $\pi$ of 
$M_I({\Bbb A})$, choose $$\eqalign{\phi\in A&(U_I({\Bbb A})M_I(F)\backslash 
G({\Bbb A}))_\pi\cap 
L^2(U_I({\Bbb A})M_I(F)\backslash G({\Bbb A}))\cr
:=&A^2(U_I({\Bbb A})M_I(F)
\backslash G({\Bbb A}))_\pi,\cr}$$ where $L^2(U_I({\Bbb A})M_I(F)\backslash 
G({\Bbb A}))$ denotes the space of $L^2$ functions on the space 
$Z_{G({\Bbb A})}
U_I({\Bbb A})M_I(F)\backslash G({\Bbb A})$. Denote the associated Eisenstein 
series by $E(\phi,\pi)\in A(G(F)\backslash G({\Bbb A}))$.
\vskip 0.30cm
\noindent
{\bf Definition.} {\it A rank $r$ 
non-abelian $L$-function $L_{F,r}^{\leq p}(\phi,\pi)$ for the function 
field $F$ associated to an
$L^2$-automorphic form $\phi\in A^2(U_I({\Bbb A})M_I(F)\backslash 
G({\Bbb A}))_\pi$ is defined by the following integration
$$L_{F,r}^{\leq p}(\phi,\pi):=
\int_{{\Cal M}_{F,r}^{\leq p}}E(\phi,\pi)(g)\,dg,
\qquad \text {Re}\,\pi\in {\Cal C}.$$}

More generally, for any standard parabolic subgroup $P_J=U_JM_J\supset P_I$ 
(so that the partition $J$ is a refinement of $I$), we have the corresponding 
relative Eisenstein series $$E_I^J(\phi,\pi)(g):=\sum_{\delta\in P_I(F)
\backslash P_J(F)}\phi(\delta g),\qquad \forall g\in P_J(F)\backslash 
G({\Bbb A}).$$ 
It is well-known that there is  an open cone ${\Cal C}_I^J$
in $\text {Re}X_{M_I}^{P_J}$ such that for $\text {Re}\,\pi\in {\Cal C}_I^J$, 
$E_I^J(\phi,\pi)\in A(P_J(F)\backslash G({\Bbb A})).$ Here $X_{M_I}^{P_J}$ 
is defined similarly as $X_M^G$ with $G$ replaced by $P_J$. Then we have 
a well-defined relative non-abelian $L$-function
$$L_{F,r}^{P_J;\leq p}(\phi,\pi):=
\int_{{\Cal M}_{F,r}^{P_J;\leq p}}E_I^J(\phi,\pi)(g)\,dg,
\qquad \text {Re}\pi\in {\Cal C}_I^J.$$ 
 \vskip 0.30cm
\noindent
{\it Remarks.} (1) Here when defining non-abelian $L$-functions we assume that 
$\phi$ comes from a single irreducible automorphic representations, but this 
restriction is rather artificial and can be removed easily: We add such a 
restriction only for the purpose of giving the 
constructions and results in a very neat way. 

\noindent
(2) We point out that the following discussion for non-abelian 
$L$-functions holds  for  relative non-abelian $L$-functions as well, 
with certain
simple modifications in a well-known manner. 
\vskip 0.30cm
\centerline
{\bf 3.4. Meromorphic Extension, Rationality and Functional Equations}
\vskip 0.30cm
With the same notation as above, set $\frak P=[\pi]$. For $w\in W$, the Weyl 
group of $G=GL_r$, fix once and for all representative $w\in G(F)$ of $w$. Set 
$M':=wMw^{-1}$ and denote the associated parabolic subgroup by $P'=U'M'$. 
$W$ acts naturally on automorphic representations, from which we obtain 
an equivalence
classes $w{\frak P}$ of automorphic representations of $M'({\Bbb A})$. As 
usual, define the associated intertwining operator $M(w,\pi)$ by 
$$(M(w,\pi)\phi)(g):=\int_{U'(F)\cap wU(F)w^{-1}\backslash U'({\Bbb A})}
\phi(w^{-1}n'g)dn',\qquad
\forall g\in G({\Bbb A}).$$ One checks that if $\langle \text {Re}\pi,
\alpha^\vee\rangle\gg 0,\forall \alpha\in\Delta_P^G$,

\noindent
(i) for a fixed $\phi$, $M(w,\pi)\phi$ depends only on the double coset 
$M'(F)wM(F)$. So $M(w,\pi)\phi$ is well-defined for $w\in W$;

\noindent
(ii) the above integral converges absolutely and uniformly for $g$ varying 
in a compact subset of $G({\Bbb A})$;

\noindent
(iii) $M(w,\pi)\phi\in A(U'({\Bbb A})M'(F)\backslash G({\Bbb A}))_{w\pi}$;
and if $\phi$ is $L^2$, which from now on we always assume, so is  
$M(w,\pi)\phi$.
\vskip 0.30cm
\noindent
{\li Basic Facts of Non-Abelian $L$-Functions.} {\it With the same notation 
above,} 

\noindent
(I) ({\bf Meromorphic Continuation}) {\it $L_{F,r}^{\leq p}(\phi,\pi)$ for 
$\text{Re}\,\pi\in {\Cal C}$ is well-defined and admits a unique 
meromorphic continuation to the whole space $\frak P$;}

\noindent
(II) ({\bf  Rationality}) {\it $L_{F,r}^{\leq p}(\phi,\pi)$ for 
$\text{Re}\,\pi\in {\Cal C}$ is a rational function on $\frak P$;}

\noindent
(III) ({\bf Functional Equations}) {\it As meromorphic functions on $\frak P$,
$$L_{F,r}^{\leq p}(\phi,\pi)=L_{F,r}^{\leq p}(M(w,\pi)\phi,w\pi),\qquad 
\forall w\in W.$$}

\noindent
{\it Proof.} This is a direct consequence of the fundamental results of 
Langlands and Morris on Eisenstein series and spectrum decompositions. 
(See e.g, [Mor1,2], [La] and/or [MW]).
Indeed, if $\phi$ is cuspidal, by definition, (I) is a direct 
consequence of Prop. II.15, Thm. IV.1.8 of [MW], 
(II) is a direct consequence of Thm. IV.1.11 of
[MW], (II) is a direct consequence of the proof of Thm. IV.1.12 of [MW]
and (II) is a direct consequence of Thm. IV.1.10 of [MW].

More generally, if $\phi$ is only $L^2$, then by Langlands and Morris' 
theory of  Eisenstein series and spectral decomposition, $\phi$ may be 
obtained as 
the residue of
relative Eisenstein series coming from cusp forms, since $\phi$ is $L^2$ 
automorphic. As such then (I), (II) and (II) are direct consequences of the 
proof  of VI.2.1(i) at p.264 of [MW].
\vskip 0.30cm
\centerline
{\bf 3.5. Holomorphicity and Singularities}
\vskip 0.30cm
Let $\pi\in \frak P=[\pi]$ and $\alpha\in\Delta_M^G$ and 
$\alpha\in R^+(T_M,G)$. Denote by $n(\alpha)$ the
smallest integer $n>0$ such that $\alpha^{*n\lambda}=1$ for all 
$\lambda\in \text {Fix}_{X_M^G}(\frak P):=\{
\nu\in X_M^G:\pi\otimes  \nu=\nu\}$ with $\alpha^*$ as defined at 
p.16-17 of [MW].
 Define then the function 
$h:\frak P\to {\Bbb C}$ by $\pi\otimes\lambda\mapsto 
\alpha^{*n(\alpha)\lambda}-1$
for all $\lambda\in X_M^G\simeq \frak a_M^G$. Set 
$H:=\{\pi'\in\frak P:h(\pi')=0\}$ and call it a root hyperplane. Clearly the 
function $h$ is determined by $H$, hence we also denote  $h$ by $h_H$. 
Note also that root hyperplanes depend on the base point $\pi$ we choose.

Let $D$ be a set of root hyperplanes. Then 

\noindent
(i) the singularities of a meromorphic function $f$ on $\frak P$ is said to 
be carried out by $D$ if for all $\pi\in\frak P$, there exist $n_\pi:D\to 
{\Bbb Z}_{\geq 0}$ zero almost everywhere such that $\pi'\mapsto 
\big(\Pi_{H\in D}h_H(\pi')^{n_\pi(H)}\big)\cdot f(\pi')$ is holomorphic at
$\pi'$;

\noindent
(ii) the singularities of $f$ are said to be without multiplicity at $\pi$ if 
$n_\pi\in\{0,1\}$;

\noindent
(iii) $D$ is said to be locally finite,  if for any compact subset 
$C\subset\frak P$, $\{H\in D:H\cap C\not=\emptyset\}$ is finite.
\vskip 0.30cm
\noindent
{\li  Basic Facts of Non-Abelian $L$-Functions.} {\it With the same notation 
above,} 

\noindent
(IV) ({\bf Holomorphicity}) (i) {\it When $\text{Re}\,\pi\in {\Cal C}$, 
$L_{F,r}^{\leq p}(\phi,\pi)$ is holomorphic;}

\noindent
(ii){\it  $L_{F,r}^{\leq p}(\phi,\pi)$ is holomorphic at $\pi$ where 
$\text{Re}\,\pi=0$;}

\noindent
(V) ({\bf Singularities}) (i) {\it There is a locally finite set of root 
hyperplanes 
$D$ such that the singularities of $L_{F,r}^{\leq p}(\phi,\pi)$ are carried 
out by $D$;}

\noindent
(ii) {\it The singularities of $L_{F,r}^{\leq p}(\phi,\pi)$ are without 
multiplicities at $\pi$ if $\langle \text{Re}\,\pi,\alpha^\vee\rangle\geq 0,
\forall \alpha\in\Delta_M^G$;}

\noindent
(iii) {\it There are only finitely many of singular hyperplanes of 
$L_{F,r}^{\leq p}(\phi,\pi)$ which intersect $\{\pi\in\frak P:
\langle\text{Re}\,\pi,\alpha^\vee\rangle\geq 0,\forall\alpha\in\Delta_M\}$.}
\vskip 0.30cm
\noindent
{\it Proof.} As above, this is a direct consequence of the fundamental 
results of 
Langlands and Morris on Eisenstein series and spectrum decompositions. 
(See e.g, [Mor1,2],  [La] and/or [MW]).
Indeed, if $\phi$ is a cusp form, (IV.i) is a direct consequence of 
Lemma IV.1.7 of [MW], while  (IV.ii)  and (IV) are direct consequence of
Prop. IV.1.11 of [MW].

In general when $\phi$ is only $L^2$ automorphic, then we have to use the 
theory of Langlands and Morris to realize $\phi$ as the residue of relative 
Eisenstein series defined using cusp forms. (See e.g., item (5) at p.198 and 
the second 
half part of p.232 of [MW].) 

As such, (IV) and (V) are direct consequence of the definition of 
residue datum and the compatibility between residue and Eisenstein series 
as stated for example under item (3) at p.263 of [MW]. 
\vskip 0.45cm
\noindent
\centerline {\li Chapter II.4. A Closed Formula for the Abelian Part}
\vskip 0.30cm
\centerline
{\bf 4.1. Modified Analytic Truncation}
\vskip 0.30cm
Let $G=GL_r$ and $P_0=M_0U_0$ be the minimal parabolic subgroup corresponding 
to the partition $(1,\cdots,1)$. Let $P_1=M_1U_1$ be a fixed standard parabolic
 subgroup with $M_1$ the standard Levi and $U_1$ the unipotent radical.

For a function field $F$ with ${\Bbb A}$ the ring of adeles, let $\pi$ be an 
irreducible automorphic representation of $M_1({\Bbb A})$. Denote by
 $A^2(U_1({\Bbb A})M_1(F)\backslash G({\Bbb A}))_\pi$ the space of 
 $L^2$-automorphic forms in the isotypic component
  $A(U_1({\Bbb A})M_1(F)\backslash G({\Bbb A}))_\pi$.
  
Then for a fixed convex polygon $p:[0,r]\to {\Bbb Q}$ and any $L^2$-automorphic
form $\phi\in A^2(U_1({\Bbb A})M_1(F)\backslash G({\Bbb A}))_\pi$ we have the 
associated non-abelian $L$-function
$$L_{F,r}^{\leq p}(\phi;\pi):=\int_{{\Cal M}_{F,r}^{\leq p}}E(\phi,\pi)(g)
\cdot d\mu(g),
\qquad \text {Re}\,\pi\in {\Cal C}$$ where $E(\phi,\pi)$ denotes the 
Eisenstein 
series associated to $\phi$ and ${\Cal C}\subset X_{M_1}^G$ is a certain 
positive cone in 3.3 over which Eisenstein series 
$E(\phi,\pi)$ converges. Recall that in 3.4, we showed that $L_{F,r}^{\leq p}
(\phi;\pi)$ admits a meromorphic continuation to the whole space 
$\frak P:=[\pi]$, the $X_{M_1}^G$ homogeneous space consisting of automorphic 
representations equivalent to $\pi$ whose typical element is 
$\pi\otimes\lambda$ with $\lambda\in X_{M_1}^G$.

On the other hand, for a suitably regular $T\in \text{Re}\frak a_M^*$, 
following Arthur, 
(see  [Ar1] and [OW],) we have the analytic truncation $\Lambda^Tf$ for 
any continuous 
function $f$ on $Z_{G({\Bbb A})}G(F)\backslash G({\Bbb A})$ defined by
$$(\Lambda^T f)(g):=\sum_P(-1)^{\text{dim}(A_P/Z_G)}
\sum_{\delta\in P(F)\backslash G(F)}f_P(\delta g)\cdot
\hat\tau_P(\log_Mm_P(\delta g)-T).$$ 
(For unknown notation, which are commonly 
used in Arthur's theory, please see [Ar1,2] and [OW].) 
Apply this 
analytic truncation to the constant function ${\bold 1}$, by Prop 1.1 
of [Ar1], we obtain a characteristic function for a certain compact subset 
in $Z_{G({\Bbb A})}G(F)\backslash G({\Bbb A})$, which we denote by 
$\Lambda^T\big(Z_{G({\Bbb A})}G(F)\backslash G({\Bbb A})\big).$
Thus, for 
$\phi\in A^2(M(F)U({\Bbb A})\backslash G({\Bbb A}))$, we have a well-defined 
integration
$$L_{F,r}^T(\phi,\pi):=
\int_{\Lambda^T\big(Z_{G({\Bbb A})}G(F)\backslash G({\Bbb A})\big)}E(\phi,\pi)
(g)
\cdot dg,\qquad \text{Re}\pi\in {\Cal C}.$$ 
Moreover, it is well-known that for analytic truncations, 
$$\Lambda^T\circ\Lambda^T=\Lambda^T$$ based on the following miracle --
By Lemma 1.1 of [Ar2], the constant term of 
$\Lambda^T\phi(x)$ along with any standard parabolic subgroup
$P_1$ is zero unless $\varpi(H_0(x)-T)<0$ for all $\varpi\in\hat\Delta_1$. 
As a direct consequence,
$$\eqalign{&L_{F,r}^T(\phi,\pi)\cr
=&
\int_{Z_{G({\Bbb A})}G(F)\backslash G({\Bbb A})}\Lambda^T{\bold 1}(g)\cdot 
E(\phi,\pi)(g)\cdot dg\cr
=&
\int_{Z_{G({\Bbb A})}G(F)\backslash G({\Bbb A})}(\Lambda^T\circ \Lambda^T)
{\bold 1}(g)\cdot E(\phi,\pi)(g)\cdot dg\cr
=&
\int_{Z_{G({\Bbb A})}G(F)\backslash G({\Bbb A})}\Lambda^T{\bold 1}(g)\cdot 
\Lambda^T E(\phi,\pi)(g)\cdot dg\cr}$$ since $\Lambda^T$ is self-adjoint.
But this latest integration is simply 
$$\int_{Z_{G({\Bbb A})}G(F)\backslash G({\Bbb A})}{\bold 1}(g)\cdot 
(\Lambda^T\circ \Lambda^T)E(\phi,\pi)(g)\cdot dg$$ since $\Lambda^T 
E(\phi,\pi)$ is rapidly decreasing and ${\bold 1}$ is of moderate growth. 
That is to say, $$L_{F,r}^T(\phi,\pi)=
\int_{Z_{G({\Bbb A})}G(F)\backslash G({\Bbb A})}\Lambda^T E(\phi,\pi)(g)\cdot 
dg.$$ 

One may try to apply such a discussion to geometric truncations as well.
For this,
attach to a fixed concave polygon $p:[0,r]\to \Bbb R$ with the property 
$p(0)=p(r)=0$ an element $T_p=(t_p^1,\cdots,t_p^r)\in \frak a_0$ by the 
conditions 
$$\lambda_i(T_p)=t_p^i-t_p^{i+1}:=[p(i)-p(i-1)]-[p(i+1)-p(i)]>0, i=1,2,
\cdots, r-1.$$
Here as usual $\{\lambda_i=e_i-e_{i+1}\}_{i=1}^{r-1}$ denotes the collection 
of positive roots of $GL_r$. Then one checks (see [We2] for details) that

(i) $T_p$ is in the positive cone of $\frak a_0$; and

(ii) $\tau_P(-H(g)-T_P)=1\Leftrightarrow p_P^g>_Pp.$

Note in particular that in (ii), $\tau_P$ instead of $\hat\tau_P$ is used. 
In other words, positive chambers rather than positive cones are used in 
geometric
truncation. We should also point out that this discussion is motivated by
 Lafforgue [Laf].
 
Moreover, following Lafforgue [Laf],  introduce a modified truncation with 
respect to a polygon $p$ by
$$(\Lambda_pf)(g):=\sum_P(-1)^{\text{dim}(A_P/Z_G)}\sum_{\delta\in 
P(F)\backslash G(F)}f_P(\delta\,g)\cdot{\bold 1}(\bar p^{\delta g}>_Pp).$$ 
Denote thus obtained moduli space (from $\Lambda_p{\bold 1}$) 
by $\Lambda_p(Z_{G({\Bbb A})}G(F)\backslash G({\Bbb A}))$. Then essentially,
the compact space $\Lambda_p(Z_{G({\Bbb A})}G(F)\backslash G({\Bbb A}))$ is
our moduli space ${\Cal M}_{F,r}^{\leq p}$ by Prop. II.2 and (ii) above. In this way, 
our problem becomes to study 
$${\Cal L}_{F,r}^{\leq p}(\phi;\pi):=
\int_{\Lambda_p(Z_{G({\Bbb A})}F(F)\backslash G({\Bbb A}))}
E(\phi,\pi)\cdot d\mu(g),\qquad\text{Re}\pi\in {\Cal C}.$$ 
\vskip 0.30cm
\centerline
{\bf  4.2. A Close Formula When $\phi$ is a Cusp Form}
\vskip 0.30cm
For general $\phi$, this turns to be a very challenging problem.
Our aim here is to see what happens for  
${L}_{F,r}^{\leq p}(\phi;\pi)$ when $\phi$ is a cusp form. 
Motivated by
the result of Langlands-Arthur on the inner product
of truncated Eisenstein series [Ar1,2] (see also [OW]), we go as follows: 

We begin with a formula for the truncated Eisenstein series. This then 
leads to the 
consideration of constant terms of Eisenstein series. While it is difficult to 
precisely describe  constant terms of Eisenstein series $E(\phi,\pi)$ 
associated with general automorphic form $\phi$, it becomes rather easy  
when $\phi$ is cuspidal. Indeed, 
for $\phi\in A_0(U_1({\Bbb A})M_1(F)\backslash G({\Bbb A}))$
and a fixed standard parabolic subgroup $P=MU$, it is well-known that
 $$E_P(\phi,\pi)(g)=\sum_{w\in W(M_1,M)}\sum_{m\in M(F)\cap wP_1(F) w^{-1}
 \backslash M(F)}\big(M(w,\pi)\phi(\pi)\big)(mg),$$ where $W(M_1,M)$ 
consisting of element  $w\in W$ such that $wM_1w^{-1}$ is a standard Levi of 
$M$ 
and $w^{-1}(\beta)>0$ for all $\beta\in R^+(T_0,M)$ and $R^+(T_0,M)$ 
denotes the set of positive roots related to $(T_0,M)$.
  
Therefore, 
$$\eqalign{&\Lambda_pE(\phi,\pi)\cr
=&\sum_{P}(-1)^{\text{dim}A_P/Z_G}
\sum_{\delta\in P(F)\backslash G(F)}E_P(\phi,\pi)(\delta g)\cdot 
{\bold 1}(\bar p^{\delta g}>_Pp)\cr
=&\sum_P(-1)^{\text{dim}A_P/Z_G}\sum_{\delta\in P(F)\backslash G(F)}\cr
&\qquad\sum_{w\in W(M_1,M)}\sum_{\xi\in M(F)\cap wP_1(F)w^{-1}\backslash M(F)}
\big(M(w,\pi)\phi\big)(\xi\delta g)\cdot {\bold 1}(\bar p^{\delta g}>_Pp).\cr}
$$
Now for any standard parabolic subgroup $P_2$, set $W(\frak a_1,\frak a_2)$ 
to be the set of 
distinct isomorphisms from $\frak a_1$ onto
$\frak a_2$ obtained by restricting elements in $W$ to $\frak a_1$, where 
$\frak a_i$ 
denotes $\frak a_{P_i}, i=1,2$
Then one checks by definition 
easily that $W(M_1,M)$ is a union over all $P_2$ of elements 
$w\in W(\frak a_1,\frak a_2)$ 
such that (i) $w\frak a_1=\frak a_2\supset \frak a_P$; and
(ii) $w^{-1}(\alpha)>0, \forall \alpha\in\Delta_{P_2}^P$.

\noindent
Hence, $$\eqalign{&\Lambda_pE(\phi,\pi)\cr
=&\sum_{P_2}\sum_{w\in W(\frak a_1,\frak a_2),P\supset P_2,w^{-1}(\alpha)>0,
\forall \alpha\in \Delta_{P_2}^P}
(-1)^{\text{dim}A_P/Z_G}\cr
&\qquad\sum_{\delta\in P(F)\backslash G(F)} {\bold 1}
(\bar p^{\delta g}>_Pp)\cdot
\sum_{\xi\in M(F)\cap wP_1(F)w^{-1}\backslash M(F)}
\big(M(w,\pi)\phi\big)(\xi\delta g)\cr
=&\sum_{P_2}\sum_{w\in W(\frak a_1,\frak a_2)}(-1)^{\text{dim}A_{P_w}/Z_G}
\sum_{\{P:P_2\subset P \subset P_w,w^{-1}(\alpha)>0,\forall \alpha\in 
\Delta_{P_2}^P\}}
(-1)^{\text{dim}A_P/A_{P_w}}\cr
&\qquad\sum_{\delta\in P(F)\backslash G(F)} {\bold 1}
(\bar p^{\delta g}>_Pp)\cdot
\sum_{\xi\in M(F)\cap wP_1(F)w^{-1}\backslash M(F)}
\big(M(w,\pi)\phi\big)(\xi\delta g).\cr}$$
where for a given $w$, we define $P_w\supset P$ by the conidition that
$$\Delta_{P_2}^{P_w}=\{\alpha\in\Delta_{P_2}:(w\pi)(\alpha^\vee)>0\}.$$
Therefore, since  $${\bold 1}
(\bar p^{\xi\delta g}>_Pp)={\bold 1}
(\bar p^{\delta g}>_Pp),\qquad \forall \delta\in P(F)\backslash G(F),
\xi\in P_2(F)\backslash P(F),$$ we have
$$\eqalign{&\Lambda_pE(\phi,\pi)\cr
=&\sum_{P_2}\sum_{w\in W(\frak a_1,\frak a_2)}(-1)^{\text{dim}A_{P_w}/Z_G}
\sum_{\{P:P_2\subset P \subset P_w,w^{-1}(\alpha)>0,
\forall \alpha\in \Delta_{P_2}^P\}}
(-1)^{\text{dim}A_P/A_{P_w}}\cr
&\qquad\sum_{\delta\in P(F)\backslash G(F)}
\sum_{\delta\in P_2(F)\backslash P(F)} \Big({\bold 1}
(\bar p^{\xi\delta g}>_Pp)\cdot
\big(M(w,\pi)\phi\big)(\xi\delta g)\Big)\cr
=&\sum_{P_2}\sum_{w\in W(\frak a_1,\frak a_2)}
(-1)^{\text{dim}A_{P_w}/Z_G}
\sum_{\{P:P_2\subset P \subset P_w,w^{-1}(\alpha)>0,
\forall \alpha\in \Delta_{P_2}^P\}}
(-1)^{\text{dim}A_P/A_{P_w}}\cr
&\qquad\sum_{\delta\in P_2(F)\backslash G(F)} \Big({\bold 1}
(\bar p^{\delta g}>_Pp)\cdot
\big(M(w,\pi)\phi\big)(\delta g)\Big)\cr
=&\sum_{P_2}\sum_{\delta\in P_2(F)\backslash G(F)}
\sum_{w\in W(\frak a_1,\frak a_2)}(-1)^{\text{dim}A_{P_w}/Z_G}
\big(M(w,\pi)\phi\big)(\xi\delta g)\cr
&\qquad\sum_{\{P:P_2\subset P \subset P_w,w^{-1}(\alpha)>0,
\forall \alpha\in \Delta_{P_2}^P\}}
(-1)^{\text{dim}A_P/A_{P_w}}
 {\bold 1}
(\bar p^{\delta g}>_Pp).\cr}$$
Set now $${\bold 1}(P_2;p;w):=
\sum_{\{P:P_2\subset P \subset P_w,w^{-1}(\alpha)>0,
\forall \alpha\in \Delta_{P_2}^P\}}
(-1)^{\text{dim}A_P/Z_G}
 {\bold 1}
(\bar p^{\delta g}>_Pp).$$
Then, we obtain the following 
\vskip 0.30cm
\noindent
{\bf Lemma.} {\it With the same notation as above,
$$\Lambda_pE(\phi,\pi)(g)=\sum_{P=MU}\sum_{\delta\in P(F)\backslash G(F)}
\sum_{w\in W(M_1), wM_1w^{-1}=M}\big(M(w,\pi)\phi\big)(\delta g)\cdot 
{\bold 1}(P;p;w)(\delta g).$$}
 
In the following calculation, I will pay no attention to the convergence: 
One may 
justify our discussion using either the standard method in [Ar1], 
and/or [OW], to first create a rapid decreasing function via 
pseudo-Eisenstein series 
or the same wave packets, then apply the inversion formula, or 
regularized integrations 
in [JLR]. Also if $\Lambda_p$ were idempotent, we would have had no chance to 
get an essential non-abelian part in our non-abelian $L$-function.

With these comments in mind, now we introduce what we call the abelian part 
${L}_{F,r}^{\leq p,\text{ab}}$ 
of our non-abelian $L$ function ${L}_{F,r}^{\leq p}$
by setting $${L}_{F,r}^{\leq p,\text{ab}}(\phi,\pi):=
\int_{Z_{G({\Bbb A})}G(F)\backslash G({\Bbb A})}\Lambda_pE(\phi,\pi)(g)\,
d\mu(g).$$ If $\Lambda_p$ were idempotent, we would have had no chance to 
get an essential non-abelian part in our non-abelian $L$-function.
It is this abelian part which we are going to calculate.

At it stands, $$\eqalign{&{L}_{F,r}^{\leq p,\text{ab}}(\phi,\pi)\cr
=&\int_{Z_{G({\Bbb A})}G(F)\backslash
 G({\Bbb A})}\sum_{P=MU}\sum_{\delta\in P(F)\backslash G(F)}
\sum_{w\in W(M_1), wM_1w^{-1}=M}\big(M(w,\pi)\phi\big)(\delta g)\cdot 
{\bold 1}(P;p;w)(\delta g)\,dg.\cr}$$ 
From an un-folding trick, it is simply
$$\eqalign{&\sum_{P}\sum_{w\in W(M_1), wM_1w^{-1}=M}
\int_{Z_{G({\Bbb A})}P(F)\backslash G({\Bbb A})}
\Big({\bold 1}(P_2;p;w)(g)\cdot
\big(M(w,\pi)\phi\big)(g)\Big)\,dg\cr
=&\sum_{P}\sum_{w\in W(M_1), wM_1w^{-1}=M}
\int_{Z_{G({\Bbb A})}U({\Bbb A})M(F)\backslash G({\Bbb A})}
\Big({\bold 1}_{P}(P;p;w)(g)\cdot
\big(M(w,\pi)\phi\big)(g)\Big)\,dg\cr}$$ where
as usual
${\bold 1}_{P}(P;p;w)(g):=\int_{U(F)\backslash U({\Bbb A})}{\bold 1}
(P;p;w)(n g)\,dn$ denotes the constant term of ${\bold 1}(P;p;w)(g)$ along $P$.

To evaluate this latest integral, we decompose it into a double
integrations over $$\eqalign{&\big(Z_{G({\Bbb A})}(Z_{M}(F)\cap Z_{M({\Bbb A})}
\backslash Z_{G({\Bbb A})}\cdot Z_{M({\Bbb A})}^1\big)\times
\big(Z_{G({\Bbb A})} Z_{M({\Bbb A})}^1U({\Bbb A)}M(F)\backslash 
G({\Bbb A})\big)\cr
=&\big(Z_{G({\Bbb A})}\cdot Z_{M({\Bbb A})}^1\backslash 
Z_{M({\Bbb A})}\big)\times
\big(Z_{M({\Bbb A})}U({\Bbb A)}M(F)\backslash G({\Bbb A})\big),\cr}$$ where
$Z_{M({\Bbb A})}^1=Z_{M({\Bbb A})}\cap M({\Bbb A})^1$.
That is to say,
$$\eqalign{{L}_{F,r}^{\leq p,\text{ab}}(\phi,\pi)
=&\sum_{P=MU}\sum_{w\in W(M_1), wM_1w^{-1}=M}
\int_{Z_{M({\Bbb A})}U({\Bbb A)}M(F)
\backslash G({\Bbb A})}dg\cr
&\cdot\int_{Z_{G({\Bbb A})}\cdot Z_{M({\Bbb A})}^1
\backslash Z_{M({\Bbb A})}} \Big({\bold 1}_{P}(P;p;w)(zg)\cdot
\big(M(w,\pi)\phi\big)(zg)\Big)dz.\cr}$$

Note now that since $X_{M_1}^G$ has no torsion, there exists a unique element 
$\pi_0$ of $\frak P:=[\pi]$ whose restriction to $A_{M_1({\Bbb A})}^G$ is 
trivial. This then allows to canonically identified $X_{M_1}^G$ with
  $\frak P$ via $\lambda_\pi\in X_{M_1}^G\mapsto \pi:=\pi_0\otimes
  \lambda_\pi\in \frak P$. Hence without loss of generality, we may simply 
  assume that the  restriction of $\pi$ to $A_{M_1({\Bbb A})}^G$ is trivial. 

Therefore, $$\eqalign{{L}_{F,r}^{\leq p,\text{ab}}(\phi,\pi)
=&\sum_{P=MU}\sum_{w\in W(M_1), 
wM_1w^{-1}=M}\int_{Z_{M({\Bbb A})}U({\Bbb A)}M(F)
\backslash G({\Bbb A})}\big(M(w,\pi)\phi\big)(g)\,dg\cr
&\qquad\cdot\int_{Z_{G({\Bbb A})}
\cdot Z_{M({\Bbb A})}^1\backslash Z_{M({\Bbb A})}} 
\Big({\bold 1}_{P}(P;p;w)(zg)
\cdot
\Big)z^{\rho_P+w\pi}dz.\cr}$$ However as $g$ may be chosen in $G({\Bbb A})^1$, 
clearly,
the integration $$\int_{Z_{G({\Bbb A})}\cdot Z_{M({\Bbb A})}^1\backslash
 Z_{M({\Bbb A})}} \Big({\bold 1}_{P}(P;p;w)(zg)\cdot
\Big)z^{\rho_P+w\pi}dz$$ is independent of $g$. Denote it by $W(P;p;w;\pi)$. 
As a direct consequence, we obtain the following
\vskip 0.30cm
\noindent
{\bf  A Closed Formula}.
{\it  With the same notation as above, for $\phi\in
 A_0(U_1({\Bbb A})M_1(F)\backslash G({\Bbb A}))_\pi$,
$${L}_{F,r}^{\leq p,\text{ab}}(\phi;\pi)=\sum_{P=MU}\sum_{w\in W(M_1), 
wM_1w^{-1}=M}
\Big(W(P;p;w;\pi)\cdot\langle M(w,\pi)\phi,1\rangle\Big).$$}
\vskip 0.30cm
\noindent
{\bf Acknowledgement.} This work is partially supported by JSPS. We would 
like to thank Deninger, Fesenko, Ueno and Zagier for their discussions, 
encouragement and interests.
\vskip 1.5cm
\noindent
\centerline {\li REFERENCES} 
\vskip 0.20cm
\noindent
[Ar1]  Arthur, J. A trace formula for reductive groups. I. Terms associated 
to classes 
in $G({\Bbb Q})$. Duke Math. J. 45  (1978), no. 4, 911--952
\vskip 0.20cm
\noindent
[Ar2] Arthur, J. A trace formula for reductive groups. II. Applications of a 
truncation 
operator. Compositio Math. 40 (1980), no. 1, 87--121.
\vskip 0.20cm
\noindent
[A] Artin, E. Quadratische K\"orper im Gebiete der h\"oheren
Kongruenzen, I,II, {\it Math. Zeit}, {\bf 19} 153-246 (1924) (See also
{\it Collected Papers}, pp. 1-94,  Addison-Wesley 1965)
\vskip 0.20cm
\noindent
[B-PGN]  Brambila-Paz, L,  Grzegorczyk, I. \&  Newstead, P.E. Geography of 
Brill-Noether Loci
for Small slops, J. Alg. Geo. {\bf 6}  645-669 (1997)
\vskip 0.20cm
\noindent
[DR] Desale, U.V. \&  Ramanan, S. Poincar\'e polynomials of the variety of
stable bundles, Math. Ann {\bf 216}, 233-244 (1975)
\vskip 0.20cm
\noindent
[HN]  Harder, G. \&  Narasimhan, M.S. On the cohomology groups of moduli spaces
of vector bundles over curves, Math Ann. {\bf 212}, (1975) 215-248
\vskip 0.20cm 
\noindent
[H]  Hasse, H. {\it Mathematische Abhandlungen}, Walter
de Gruyter, Berlin-New York, 1975.
\vskip 0.20cm 
\noindent
[JLR] Jacquet, H.; Lapid, E.\& Rogawski, J.  Periods of automorphic forms. 
J. Amer. Math. Soc. 12 (1999), no. 1, 173--240.
\vskip 0.20cm
\noindent
[Laf] Lafforgue, L. {\it Chtoucas de Drinfeld et conjecture de 
Ramanujan-Petersson}. 
Asterisque No. 243 (1997)
\vskip 0.20cm 
\noindent
[LS] Lagarias, J. \& Suzuki, in preparation
\vskip 0.20cm
\noindent
[La] Langlands, R. {\it On the functional equations satisfied by Eisenstein 
series}, 
Springer LNM {\bf 544}, 1976
\vskip 0.20cm
\noindent
[MW] Moeglin, C. \& Waldspurger, J.-L. {\it Spectral decomposition and 
Eisenstein series}. 
Cambridge Tracts in Mathematics, {\bf 113}. Cambridge University Press, 
Cambridge, 1995.
\vskip 0.20cm
\noindent 
[Mo] Moreno, C. {\it  Algebraic curves over finite fields.}
Cambridge Tracts in Mathematics, {\bf 97}, Cambridge University Press, 1991
\vskip 0.20cm
\noindent 
[Mor1] Morris, L.E. Eisenstein series for reductive groups over global 
function fields I: the cusp form case, Can. J. Math. {\bf 34}, (1982), 91-168
\vskip 0.20cm
\noindent 
[Mor2] Morris, L.E. Eisenstein series for reductive groups over global 
function fields II: the general case, Can. J. Math. {\bf 34}, (1982), 1112-1182
\vskip 0.20cm
\noindent 
[Mu] Mumford, D. {\it Geometric Invariant Theory}, Springer-Verlag,
(1965)
\vskip 0.20cm
\noindent
[NR] Narasimhan, M.S. \&  Ramanan, S. Moduli of vector bundles on a compact 
Riemann surfaces,
Ann. of Math. {\bf 89} 14-51 (1969)
\vskip 0.20cm
\noindent
[NS] Narasimhan, M. S. \& Seshadri, C. S. Stable and unitary vector bundles 
on a compact Riemann surface. 
Ann. of Math. (2) 82 1965 540--567. 
\vskip 0.20cm 
\noindent
[OW] Osborne, M. \& Warner, G. The Selberg trace formula. II. Partition, 
reduction, 
truncation. Pacific J. Math. 106 (1983), no. 2, 307--496.
\vskip 0.20cm
\noindent
[Se]  Seshadri, C.S. {\it Fibr\'es vectoriels sur les courbes alg\'ebriques},
 Asterisque {\bf
96}, 1982
\vskip 0.20cm
\noindent 
[Ta] Tate, J. Fourier analysis in number fields and Hecke's
zeta functions, Thesis, Princeton University, 1950 
\vskip 0.20cm
\noindent 
[W1] Weil, A. {\it Sur les courbes alg\'ebriques et les vari\'et\'es qui
s'en d\'eduisent}, Herman, Paris (1948)
\vskip 0.20cm
\noindent
[W2] Weil, A. {\it Adeles and algebraic groups}, Prog. in Math, {\bf 23} (1982)
\vskip 0.20cm
\noindent
[We1] Weng, L.  A Program for Geometric Arithmetic, at http://xxx.lanl.gov/abs/math.AG/0111241
\vskip 0.20cm
\noindent 
[We2] Weng, L. Non-Abelian $L$-Function for Number Fields, 2003
\vskip 0.20cm
\noindent
[We3] Weng, L.  Refined Brill-Noether Locus and Non-Abelian Zeta
Functions  for Elliptic  Curves, in {\it Proceedings of Algebraic Geometry 
in South-East Asia},
2002, World Sci.
\vskip 0.20cm
\noindent
[Z] Zagier, D. The Rankin-Selberg method for automorphic functions 
which are not of rapid decay, J. Fac. Sci. Univ. Tokyo 28 (1982), 415--437. 
\end